\documentclass{article}
\usepackage{amsmath}
\usepackage{amsfonts}
\usepackage{amsthm}
\usepackage{amssymb}
\usepackage{centernot}
\usepackage{url}

\numberwithin{equation}{section}

\def\lcm{{\rm lcm}}

\def\intZ{\mathbb{Z}}

\newtheorem{Theorem}{Theorem}[section] 
\newtheorem{Lemma}{Lemma}

\newtheorem*{Corollary_2}{Corollary to Lemma 2}
\newtheorem*{Corollary3_1*}{Corollary to Proposition 3.1}
\newtheorem*{Corollary4_1*}{Corollary to Proposition 4.1}
\newtheorem*{CorollaryThm4_4*}{Corollary to Theorem 4.4}
\newtheorem*{Theorem4_3}{Theorem 4.3}
\newtheorem*{Theorem4_4}{Theorem 4.4}
\newtheorem*{Theorem4_5}{Theorem 4.5}
\newtheorem*{Observation4_6}{Observation 4.6}
\newtheorem*{Theorem4_7}{Theorem 4.7}
\newtheorem*{Theorem4_8}{Theorem 4.8}

\newtheorem{Proposition}{Proposition}[section]
\newtheorem{Criterion}{Criterion}

\def\rad{{\rm rad}}

\begin{document}

\begin{center}
{\large\bf 
Number of solutions to $a^x+b^y=c^z$ with $\gcd(a,b)>1$
}
  
\bigskip 

Reese Scott  

Robert Styer

\bigskip

Keywords: {ternary purely exponential Diophantine equation, number of solutions}

2020 Subject Class: {11D61}

\end{center}


18 December 2024  

\bigskip

\begin{abstract}  
We show that there are at most two solutions in positive integers $(x,y,z)$ to the equation $a^x+b^y=c^z$ for positive integers $a$, $b$, and $c$ all greater than one, with just one exceptional case when $\gcd(a,b)=1$, and just one exceptional infinite family of cases when $\gcd(a,b)>1$ (two solutions $(x_1,y_1,z_1)$ and $(x_2,y_2,z_2)$ are considered the same solution if $\{ a^{x_1},  b^{y_1} \} = \{ a^{x_2}, b^{y_2} \}$).  The case in which $\gcd(a,b)=1$ has been handled in a series of successive results by Scott and Styer, Hu and Le, and Miyazaki and Pink, who showed that there are at most two solutions, excepting $(\{a,b\},c) = (\{3,5\},2)$, which gives three solutions.  So here we treat the case $\gcd(a,b)>1$, showing that in this case there are at most two solutions, excepting $(a,b,c) = (2^u, 2^v, 2^w)$ with $\gcd(uv,w)=1$, which gives an infinite number of solutions.  
This generalizes work of Bennett, who proved, for both $\gcd(a,b)=1$ and $\gcd(a,b)>1$, there are at most two solutions $(y,z)$ to the equation $a + b^y = c^z$, and conjectured there are exactly eleven $(a,b,c)$ giving two solutions to this equation (assuming $b$ and $c$ are not perfect powers).  

For both $\gcd(a,b)=1$ and $\gcd(a,b)>1$, there are an infinite number of $(a,b,c)$ giving two solutions $(x,y,z)$ to the title equation, which are described in detail in this and a cited previous paper. 

In a further result, in which we no longer say that two solutions $(x_1,y_1,z_1)$ and $(x_2,y_2,z_2)$ are considered the same solution if $\{ a^{x_1},  b^{y_1} \} = \{ a^{x_2}, b^{y_2} \}$, we list all cases with more than two solutions.  
\end{abstract}

MSC: 11D61

\section{Introduction}   

For integers $a$, $b$, and $c$ all greater than one, we consider $N(a,b,c)$, the number of solutions in positive integers $(x,y,z)$ to the equation
$$ a^x + b^y = c^z.  \eqno{(1.1)}$$
In this paper, we treat the case in which $\gcd(a,b)>1$, but first we give a brief history of previous results for the case in which $\gcd(a,b)=1$.  

For $(a,b,c)$ with $\gcd(a,b)=1$, an effective upper bound for $N(a,b,c)$ was first given by A. O. Gel'fond \cite{G} (Mahler \cite{M} had earlier shown that the number of solutions was finite, using his $p$-adic analogue of the Diophantine approximation method of Thue-Siegel, but his method is ineffective).  Hirata-Kohno \cite{HK} used an application of an upper bound on the number of solutions of binary $S$-unit equations due to F. Beukers and H. P. Schlickewei \cite{BSch} to obtain $N(a,b,c) \le 2^{36}$, later improved to $N(a,b,c) \le 200$ (unpublished).  The following more realistic upper bounds for $N(a,b,c)$ when $\gcd(a,b)=1$ have been obtained in recent years:

(1)  (R. Scott and R. Styer \cite{ScSt6a})  If $ 2 \nmid c$ then $N(a,b,c) \le 2$.

(2)  (Y. Z. Hu and M. H. Le \cite{HL1})  If $\max\{ a,b,c\} > 5 \cdot 10^{27}$, then $N(a,b,c) \le 3$.  

(3)   (Y. Z. Hu and M. H. Le \cite{HL2})  If $2 \mid c$ and $\max\{ a,b,c\} > 10^{62}$, then $N(a,b,c) \le 2$.

(4)  (T. Miyazaki and I. Pink \cite{MP})  If $2 \mid c$ and $\max\{a,b,c \} \le 10^{62}$, then $N(a,b,c) \le 2$ except for $N(3,5,2) = N(5,3,2) = 3$.  

Noting that (1), (3), and (4) show that (1.1) has at most two solutions $(x,y,z)$ when $\gcd(a,b)=1$ except for the case $(\{a,b\},c)=(\{3,5\},2)$, in what follows we consider $\gcd(a,b)>1$.  We know of no results on the number of solutions to (1.1) for the case $\gcd(a,b)>1$ for general $a$, $b$, $c$, except in the special case (known as the Pillai case) in which one of $x$ or $y$ is a fixed positive integer, in which case the proof that there are at most two solutions $(y,z)$ (respectively, $(x,z)$) is extremely short and straightforward (see \cite{Be}).

We will show that (1.1) has at most two solutions when $\gcd(a,b)>1$, provided $a$, $b$, and $c$ are not all powers of 2.  In counting the number of solutions to (1.1) to determine $N(a,b,c)$ we use the following:

\begin{Criterion}    
Two solutions $(x_1,y_1, z_1)$ and $(x_2,y_2,  z_2)$ are considered the same solution if $\{ a^{x_1}, b^{y_1} \}  =  \{ a^{x_2}, b^{y_2} \}$.  
\end{Criterion}

When neither $a$ nor $b$ is a perfect power, Criterion 1 becomes the following:

\smallskip
\noindent
{\it Two solutions $(x_1,y_1, z_1)$ and $(x_2,y_2,  z_2)$ are considered the same solution if $a=b$ and $\{ {x_1}, {y_1} \}  =  \{ {x_2}, {y_2} \}$.   }
\smallskip

When we say that two solutions $(x_1,y_1,z_1)$ and $(x_2,y_2,z_2)$ are {\it distinct}, we do not mean that they are counted as different solutions under Criterion 1, but rather simply that the ordered triple $(x_1,y_1,z_1)$ is not the same as the ordered triple $(x_2,y_2,z_2)$.  In general, even if it is not explicitly stated, two solutions $(x_1,y_1,z_1)$ and $(x_2,y_2,z_2)$ are assumed to be distinct in this sense.  On the other hand, if we are concerned with whether or not solutions are different under Criterion 1, we use the notation $N(a,b,c)$.   

Before proceeding, we discuss all known cases with $N(a,b,c)>1$. 

Considering $(a,b,c)$ with $\gcd(a,b)=1$ for which (1.1) has two solutions $(x_1,y_1,z_1)$ and $(x_2,y_2,z_2)$, we find one infinite family of such $(a,b,c,x_1,y_1,z_1,x_2,y_2,z_2)$ (given by $\{a,b\} = \{2, 2^{n} -1 \}$, $c=2^n+1$, $n >1$) and several other such $(a,b,c,x_1,y_1,z_1, x_2,y_2,z_2)$, all but one of which are apparently anomalous cases not in an infinite family (the case $\{a,b\} = \{2,7\}$ with $c=3$ gives a member of the known infinite family).  Only one case  $(\{ a,b\} = \{ 3,5\}, c=2)$ gives more than two solutions $(x,y,z)$.  See \cite{ScSt6a} for a list of the known cases with two solutions to (1.1) when $\gcd(a,b)=1$; it is conjectured in \cite{ScSt6a} that there are no further $(a,b,c)$ giving more than one solution to (1.1).  Although this conjecture is far from solved, it has been shown that for any $(a,b,c)$ with $\gcd(a,b)=1$ giving more than one solution to (1.1) with $(a,b,c)$ not listed in \cite{ScSt6a} we must have $c \ne 2$ \cite{Sc}, $c$ not a  Fermat prime \cite{MP2}, $c \ne 6$ \cite{MP2}, $c \ne 13$ \cite{MP3}.  If we restrict $a$, $b$, and $c$ to be primes or prime powers, then any $(a,b,c)$ with $\gcd(a,b)=1$ giving more than one solution to (1.1) with $(a,b,c)$ not listed in \cite{ScSt6a} must have $c>10^{18}$ \cite{LSS}.    

Considering $(a,b,c)$  with $\gcd(a,b)>1$ for which (1.1) has two solutions $(x_1,y_1,z_1)$ and $(x_2,y_2,z_2)$ which are considered different solutions under Criterion 1, we find four infinite families of such $(a,b,c, \allowbreak x_1,y_1,z_1, \allowbreak  x_2,y_2,z_2)$ and ten anomalous such $(a,b,c, \allowbreak x_1,y_1,z_1, \allowbreak  x_2,y_2,z_2)$ which are apparently not in an infinite family.  (In counting the anomalous $(a,b,c, \allowbreak x_1,y_1,z_1, \allowbreak  x_2,y_2,z_2)$  we assume that none of $a$, $b$, or $c$ is a perfect power, and we further assume that $a<b$, $x_1<x_2$ or, when $x_1=x_2$, $y_1<y_2$.)  The question arises:  are there any further $(a,b,c, x_1,y_1,z_1, x_2,y_2,z_2)$ with $\gcd(a,b)>1$ not included among the ten anomalous nine-tuples or the four infinite families?  We will address this question in Section 5.       

The four infinite families mentioned in the previous paragraph are: 

(i).  $(a,b,c) = (2,2^u (2^{h-1} - 1),2^u (2^{h-1} + 1))$, $(x_1,y_1,z_1) = (u+1, 1, 1)$ and $(x_2,y_2,z_2)=(2u+h+1, 2, 2)$, $u>0$, $h>1$.  

(ii).  $(a,b,c) = (2 \cdot 3^t, 3, 3)$, $(x_1,y_1,z_1) = (1, t, t+1)$ and $(x_2,y_2,z_2) = (3, 3t, 3t+2)$, $t>0$.

(iii).  $(a,b,c) = (g^j, g^{ju} d, g^{ju} (d+1))$, $(x_1,y_1,z_1) = (u,1,1)$ and $(x_2,y_2,z_2) = (ku+\frac{w}{j}, k, k)$ where $j>0$, $u>0$, $d>0$, $2 \nmid g > 1$, $(d+1)^k - d^k = g^w$, $w>0$, $k>1$.  

(iv).  $(a,b,c) = (2^i g^{j}, 2^{iu - 1} g^{ju} d, 2^{iu-1} g^{ju} (d+2))$, $(x_1,y_1,z_1) = (u, 1, 1)$ and $(x_2,y_2,z_2) = (ku+\frac{w}{j}, k, k)$, where $i>0$, $j>0$,  $u>0$, $d > 0$, $2 \nmid d$, $2 \nmid g > 1$,  $g^w$ is the greatest odd divisor of $(d +2)^k - d^k$, $w>0$, $2 \mid k > 0$, and $ k - v = h-(iw/j)$, where $2^h \parallel 2 d +2$ and $2^v \parallel k$. 

Before discussing these four infinite families, we point out that readers concerned only with showing $N(a,b,c) \le 2$ need not be concerned with the details relating to infinite families (see Comment 4.4 in Section 4). The purpose of considering infinite families is two-fold: first, to list all cases with more than two solutions $(x,y,z)$ to (1.1) when Criterion 1 is not used, and, second, to explore {\it anomalous} $(a,b,c, \allowbreak x_1,y_1,z_1, \allowbreak  x_2,y_2,z_2)$, which we define below. 

We now consider the infinite families (i), (ii), (iii), and (iv).  
Note that the infinite family (i) can be viewed as a subset of the infinite family (iv): redefine $g$, taking $g=i=j=1$, which requires $d = 2^{h-1} -1$ and $k = 2$, which gives $w = h-1$, and note that $u_{(iv)} = u_{(i)} + 1$ where $u_{(iv)} = u$ as in (iv) and $u_{(i)} = u$ as in (i).  But for the purposes of Propositions 3.1, 3.4, 4.1, and 4.2, it is clearer to treat (i) as a distinct family.  (We take $g>1$ in (iv) so that (i) is distinct from (iv).)

We give three definitions concerning the infinite families (i), (ii), (iii), (iv).  

Definition 1: For a given infinite family, define $F$ as follows: choose an allowable integer value for each of the parameters in the infinite family to create an ordered nine-tuple of fixed positive integers $(a,b,c, \allowbreak x_1,y_1,z_1, \allowbreak  x_2,y_2,z_2)$; then, for this infinite family, $F$ is the set of all such nine-tuples.  

Definition 2:  We say that the solution $(x,y,z)$ to (1.1) for the triple $(a,b,c)$ corresponds to the solution $(X, Y, Z)$ to (1.1) for the triple $(A,B,C)$ if $\{ a^x, b^y \} = \{ A^X,  B^Y \}$.  

Comment:  Note that for any $F$ as in Definition 1, $(x_1,y_1,z_1)$ does not correspond to $(x_2,y_2,z_2)$.  

Definition 3: Let $a$, $b$, $c$, $x_1$, $y_1$, $z_1$, $x_2$, $y_2$, $z_2$ be positive integers such that $(x_1,y_1,z_1)$ and $(x_2,y_2,z_2)$ are distinct solutions to (1.1) for the triple $(a,b,c)$.  We say that $(a,b,c, \allowbreak  x_1, y_1, z_1, \allowbreak  x_2, y_2, z_2)$ is in a given infinite family if there exists a nine-tuple $(A,B,C, \allowbreak X_1,Y_1,Z_1, \allowbreak X_2,Y_2,Z_2) \in F$ such that each of the solutions $(X_1,Y_1,Z_1)$, $(X_2,Y_2,Z_2)$ to (1.1) for the triple $(A,B,C)$ corresponds to one of the solutions $(x_1,y_1,z_1)$, $(x_2,y_2,z_2)$ to (1.1) for the triple $(a,b,c)$.  

We give an example to show how these three definitions are used.  Consider (1.1) with $a = b = 7$ and $c= 98$.  We find three solutions $(x,y,z)$, the last two of which correspond to each other so that they are considered the same solution under Criterion 1 (we will later show there are no further solutions): 
$$ (a,b,c) = (7,7, 98): (x_1,y_1,z_1) = (2,2,1), (x_2,y_2,z_2) = (6,7,3), (x_3,y_3,z_3) = (7,6,3). \eqno{(1.2)}$$
Now consider (1.1) with $a=7$, $b=49$, and $c=98$.  We find two solutions (we will later show there are no further solutions):
$$ (a,b,c) = (7, 49, 98): (x_1,y_1,z_1) = (2,1,1), (x_2,y_2,z_2) = (7,3,3). \eqno{(1.3)}$$

Note that $a$, $b$, $c$, $x_1$, $y_1$, $z_1$, $x_2$, $y_2$, $z_2$ in (1.3) exactly match the values given for $a$, $b$, $c$, $x_1$, $y_1$, $z_1$, $x_2$, $y_2$, $z_2$  in the infinite family (iii) when $g = 7$, $j=1$, $u=2$, $d=1$, $k=3$, and $w=1$, so that $(a,b,c, x_1, y_1, z_1, x_2, y_2,z_2) = (7,49,98,2,1,1,7,3,3)$ is in $F$ for the infinite family (iii).

On the other hand, $a$, $b$, $c$, $x_1$, $y_1$, $z_1$, $x_2$, $y_2$, $z_2$ as in (1.2) do not match the values given for $a$, $b$, $c$, $x_1$, $y_1$, $z_1$, $x_2$, $y_2$, $z_2$ in the infinite family (iii) for any choice of $g$, $j$, $u$, $d$, $k$, $w$ (due both to the choice of $b$ and the order of the $(x_2,y_2)$), so that $(a,b,c,x_1,y_1,z_1,x_2,y_2,z_2) = (7,7,98,2,2,1,6,7,3)$ is not in $F$ for the infinite family (iii).  However, we still say that $(a,b,c,x_1,y_1,z_1,x_2,y_2,z_2) = (7,7,98,2,2,1,6,7,3)$ is in the infinite family (iii) since the solutions $(2,2,1)$ and $(6,7, 3)$ for the triple $(a,b,c) = (7,7,98)$ correspond respectively to the solutions $(2,1,1)$ and $(7,3,3)$ for the triple $(a,b,c) = (7, 49, 98)$.  

We now consider nine-tuples $(a,b,c, \allowbreak x_1,y_1,z_1, \allowbreak  x_2,y_2,z_2)$ not in any of the four infinite families.  We define {\it anomalous} $(a,b,c, \allowbreak x_1,y_1,z_1, \allowbreak  x_2,y_2,z_2)$ as nine-tuples such that $\gcd(a,b) > 1$ and $(a,b,c)$ gives two solutions   $(x_1,y_1,z_1)$ and $(x_2,y_2,z_2)$ to (1.1) which do not correspond to each other, with $(a,b,c, \allowbreak x_1,y_1,z_1, \allowbreak  x_2,y_2,z_2)$ not in any of the infinite families (i), (ii), (iii), (iv).  
In Section 4 we will use infinite families to show that, for any such anomalous $(a,b,c, \allowbreak x_1,y_1,z_1, \allowbreak  x_2,y_2,z_2)$, $(a,b,c)$ gives no further solutions $(x,y,z)$ to (1.1) even if Criterion 1 is not used. 

The only known such anomalous $(a,b,c,x_1,y_1, z_1,x_2,y_2,z_2)$ are the ten anomalous nine-tuples with $\gcd(a,b)>1$ mentioned above, which are: 

$(a,b,c) = (2, 6, 38)$, $(x_1,y_1,z_1) = (1,2,1)$ and $(x_2,y_2,z_2)= (5,1,1)$.

$(a,b,c) = (3, 6, 15)$, $(x_1,y_1,z_1) = (2,1,1)$ and $(x_2,y_2,z_2)=(2,3,2)$.

$(a,b,c) = (6, 15, 231)$, $(x_1,y_1,z_1) = (1,2,1)$ and $(x_2,y_2,z_2)=(3,1,1)$.

$(a,b,c) = (3, 1215, 6)$, $(x_1,y_1,z_1) = (4,1,4)$ and $(x_2,y_2,z_2)=(8,1,5)$.

$(a,b,c) = (3, 6, 7857)$, $(x_1,y_1,z_1) = (4,5,1)$ and $(x_2,y_2,z_2)=(8,4,1)$.

$(a,b,c) = (5, 275, 280)$, $(x_1,y_1,z_1) = (1,1,1)$ and $(x_2,y_2,z_2)=(7,1,2)$.

$(a,b,c) = (5, 280, 78405)$, $(x_1,y_1,z_1) = (1,2,1)$ and $(x_2,y_2,z_2)=(7,1,1)$.

$(a,b,c) = (30,70,4930)$, $(x_1,y_1,z_1) = (1,2,1)$ and $(x_2,y_2,z_2)=(5,2,2)$.

$(a,b,c) = (30, 4930, 24304930)$, $(x_1,y_1,z_1) = (1,2,1)$ and $(x_2,y_2,z_2)=(5,1,1)$.

$(a,b,c) = (2, 88, 6)$, $(x_1,y_1,z_1) = (7,1,3) $ and $(x_2,y_2,z_2)=(5,2,5)$.

Nine of these anomalous cases can be derived from the five cases with $\gcd(a,b)>1$ listed by Bennett \cite{Be} in his list of eleven double solutions for the Pillai case.  Four of these five Pillai cases generate two items on our list (since Pillai equations can be rearranged) but one of them does not since one of its two possible arrangements is a member of the infinite family (i).  The tenth anomalous case is derived not from one of the known Pillai cases with $\gcd(a,b)>1$ but rather from the equations $1 + 2 \cdot 11^2 = 3^5$, $2^4 + 11 = 3^3$.  

We have not found any further anomalous $(a,b,c,x_1,y_1, z_1,x_2,y_2,z_2)$ with $2 \le \gcd(a,b)<1050$, $a/\gcd(a,b)<1050$, $b/\gcd(a,b)<10^5$, $a^x < 10^{50}$, $b^y<10^{50}$.  In Section 5 we show that any further anomalous $(a,b,c,x_1,y_1,z_1,x_2,y_2,z_2)$ must have $\rad(abc)>10^7$.  

Our main result, Theorem 1.1 which follows, is an immediate consequence of Theorem 4.3 in Section 4.  Theorem 4.5 gives additional information on anomalous cases.  Theorem 4.7 uses infinite families to  treat (1.1) without using Criterion 1, listing all cases with more than two solutions $(x,y,z)$.   Theorem 4.8 will give a refined version of Theorem 1.1, made possible by using results on infinite families.

\begin{Theorem}  
 If $a$, $b$, $c$ are positive integers all greater than one with at least one of $a$, $b$, $c$ not a power of 2 and with $(\{ a,b \}, c) \ne  (\{ 3,5 \}, 2)$, then $N(a,b,c) \le 2$.
\end{Theorem}

For the case $\gcd(a,b)=1$, this result was first stated by Miyazaki and Pink \cite{MP}, who handled the case $2 \mid c$ with $\max\{a,b,c\} \le 10^{62}$, which completed the treatment of the case $\gcd(a,b)=1$ since the case $2 \nmid c$ and the case $\max\{ a,b,c\} > 10^{62}$ with $ 2 \mid c$ had already been handled in \cite{ScSt6a} and \cite{HL2}, respectively. So to prove Theorem 1.1 it suffices to consider only the case $\gcd(a,b)>1$.  We will need several lemmas and preliminary propositions which follow.

\section{Lemmas}   

For given integers $a$, $b$, and $c$ all greater than one with $\gcd(a,b)>1$, let $Q$ be the set of all primes $p$ such that $p \mid a$, $p \mid b$, and $p \mid c$.  For every prime $p \in Q$ let $p^{\alpha_p} \parallel a$, $p^{\beta_p} \parallel b$,  $p^{\gamma_p} \parallel c$.  Let $a_1$ be the greatest divisor of $a$ not divisible by any prime in $Q$, and define $b_1$ and $c_1$ similarly for $b$ and $c$.

\begin{Lemma}  
Let $(x,y,z)$ be a solution to (1.1).  Then, for $p \in Q$, two members of the triple $\{ \alpha_p x, \beta_p y, \gamma_p z\}$ must be equal and these two members must be less than or equal to the third member.  
\end{Lemma}

\begin{proof} 
The lemma is an immediate consequence of the above definitions.
\end{proof}

We say that a solution to (1.1) is Type A for $p$ (where $p \in Q$) if $\alpha_p x > \beta_p y = \gamma_p z$; we say a solution to (1.1) is Type B for $p$ if $\beta_p y > \alpha_p x = \gamma_p z$; we say a solution to (1.1) is Type C for $p$ if $\gamma_p z > \alpha_p x = \beta_p y$; we say that a solution to (1.1) is Type O for $p$ if $\alpha_p x = \beta_p y = \gamma_p z$.

\begin{Lemma} 
Let $p$ and $q$ be two primes in the set $Q$ such that $\frac{\alpha_p}{\beta_p} > \frac{\alpha_q}{\beta_q}$ and 
$\frac{\alpha_p}{\gamma_p} > \frac{\alpha_q}{\gamma_q}$.  Then (1.1) can have no solutions which are Type B, C, or O for $p$.
\end{Lemma}

\begin{proof}
If (1.1) has a solution $(x,y,z)$ which is Type B, C or O for $p$ then $\alpha_p x \le \beta_p y$ and $\alpha_p x \le \gamma_p z$,  so that $\alpha_q x < \beta_q y$ and $\alpha_q x < \gamma_q z$, contradicting Lemma 1.    
\end{proof}

\begin{Corollary_2}  
If (1.1) has a solution which is Type A for some prime $p \in Q$, and if this solution is not Type A for some prime $q \in Q$, then (1.1) can have no solutions which are Type B, C, or O for $p$.  
\end{Corollary_2}

\begin{proof} 
If a solution $(x,y,z)$ to (1.1) is Type A for some $p \in Q$, and is not Type A for some prime $q \in Q$, then $\alpha_p x > \beta_p y$, $\alpha_p x > \gamma_p z$, $\alpha_q x \le \beta_q y$, and $\alpha_q x \le \gamma_q z$.  So $\frac{\alpha_p}{\beta_p} > \frac{\alpha_q}{\beta_q}$ and 
$\frac{\alpha_p}{\gamma_p} > \frac{\alpha_q}{\gamma_q}$, and the corollary follows from Lemma 2.
\end{proof}

We will also need several general elementary results, which follow.

\begin{Lemma}  
Let $R$, $S$, $M$, and $t_1$ be positive integers such that $R>S$, $\gcd(R,S)=1$, and $M \mid R^{t_1} - (-1)^\epsilon  S^{t_1}$ for a fixed choice of $\epsilon \in \{ 0, 1\}$.  Let $t_0$ be the least positive integer such that $M \mid R^{t_0} - (-1)^\epsilon S^{t_0}$ for this choice of $\epsilon$.  Then $t_0 \mid t_1$. 
\end{Lemma}

\begin{proof}
For $\epsilon = 0$ this is Lemma 3.1 of \cite{ScSt6}.  A similar method of proof handles $\epsilon = 1$: take $M > 2$ (since the lemma clearly holds for $M \le 2$) and let $t_0$ and $t_1$ be as in the formulation of the lemma; if we assume $t_0 \nmid t_1$, we can let $t_1 = s t_0 + r$ where $2 \nmid s$, $0 < r < 2 t_0$, and $r \ne t_0$, so that $R^r \equiv S^r \bmod M$; we can let $r = t_0 \pm r_1$ where $0 < r_1 < t_0$; then $R^{r_1} \equiv - S^{r_1} \bmod M$, contradicting the definition of $t_0$.  
\end{proof}

The following lemma sharpens Lemma 3.2 of \cite{ScSt6}.  

For any positive integer $m$ we define $P(m)$ to be the set of primes which divide $m$.  

\begin{Lemma}  
If $R$, $S$, $n_1$, and $n_2$ are positive integers with $\gcd(R,S)=1$, $R>S$, $n_1 < n_2$, and each prime dividing $R^{n_2} - S^{n_2} $ also divides $R^{n_1} - S^{n_1}$, then $n_1=1$, $n_2=2$ and $R+S=2^h$ for some integer $h > 1$.
\end{Lemma}

\begin{proof} 
Let $n_0$ be the least number such that $\rad(R^{n_2} - S^{n_2}) \mid R^{n_0} - S^{n_0}$.  By Lemma 3 we have $n_0 \mid n_2$, so that $R^{n_0} - S^{n_0} \mid R^{n_2} - S^{n_2}$, so that we have $P(R^{n_0} - S^{n_0}) = P(R^{n_2} - S^{n_2})$.  We will show that $n_2 = 2 n_0$.  

Assume $p$ is an odd prime which divides $n_2/n_0$.  Then $R^{n_0} -S^{n_0} \mid R^{p n_0} - S^{p n_0} \mid R^{n_2} - S^{n_2}$, so that $P(R^{p n_0} - S^{p n_0}) = P(R^{n_0} - S^{n_0})$.  Since $n_0 \mid n_1 < n_2$, we see that $n_0$ and $n_2$ are distinct.  Consider $R^{p n_0} - S^{p n_0} = ( (R^{n_0} - S^{n_0}) + S^{n_0})^{p} - S^{p n_0}$.  From the binomial expansion of $ ( (R^{n_0} - S^{n_0}) + S^{n_0})^{p}$ we see that if $p \nmid R^{n_0} - S^{n_0}$ then $( R^{p n_0} - S^{p n_0})/(R^{n_0} - S^{n_0}) > 1$ is prime to $R^{n_0} - S^{n_0}$, and, if $p \mid R^{n_0} - S^{n_0}$ then $( R^{p n_0} - S^{p n_0})/( p(R^{n_0} - S^{n_0}) ) > 1$ is prime to $R^{n_0} - S^{n_0}$; in either case we have a contradiction to $P(R^{n_0} - S^{n_0}) = P(R^{p n_0} - S^{p n_0})$.  So ${n_2}/{n_0}$ is not divisible by any odd prime.  

If $4 \mid {n_2}/{n_0}$, then $R^{2 n_0} + S^{2 n_0} \mid R^{n_2} - S^{n_2}$, again giving a contradiction since $R^{2 n_0} + S^{2 n_0}$ is divisible by an odd prime which is prime to $R^{n_0} - S^{n_0}$.  
So we must have $n_2 = 2 n_0$, so that, since $P(R^{n_0} - S^{n_0})= P(R^{n_2} - S^{n_2})$, we have $R^{n_0} + S^{n_0} = 2^h$ for some $h>1$, which requires $n_0=1$.  Since $1 = n_0 \le n_1 < n_2 = 2 n_0 = 2$, we have $n_1 = n_0$, giving Lemma 4.
\end{proof}

\begin{Lemma} 
If $R$, $S$, $n_1$, and $n_2$ are positive integers with $\gcd(R,S)=1$, $R>S$, $n_1 < n_2$, and each prime dividing $R^{n_2} + S^{n_2} $ also divides $R^{n_1} + S^{n_1}$, then $(R, S, n_1, n_2) = (2, 1, 1, 3)$. 
\end{Lemma}

\begin{proof} 
The proof of Lemma 5 follows that of Lemma 4 with $n_0$ redefined to be the least number such that $\rad(R^{n_2} + S^{n_2}) \mid R^{n_0} + S^{n_0}$: noting that $2 \nmid {n_2}/{n_0}$ and considering the  binomial expansion of 
$((R^{n_0} + S^{n_0})- S^{n_0})^p$, we see that the only possibility is  $R^{n_0} = 2$, $S^{n_0} = 1$, $n_2 = 3 n_0$.
So $1 = n_0 \le n_1 < n_2 = 3$, so that, since $2 \nmid ({n_1}/{n_0})$, we have $n_1 = n_0$.  
\end{proof}

\begin{Lemma}  
Let $R$, $S$, $n_1$, and $n_2$ be positive integers with $\gcd(R,S)=1$, $R>S$, and $n_1 \mid n_2$.  Let $p$ be a prime such that $p^{v_1} \parallel R^{n_1} - S^{n_1}$ where $p^{v_1}>2$.  If $p^{v_2} \parallel R^{n_2} - S^{n_2}$, then $p^{v_2 - v_1} \mid \frac{n_2}{n_1}$.     
\end{Lemma}

\begin{proof}  
Since $n_1 \mid n_2$, we have $v_1 \le v_2$.  
Considering the binomial expansion of $((R^{n_1} - S^{n_1}) + S^{n_1})^{k}$, we see that $k=p$ is the least value of $k$ such that $p^{v_1+1} \mid R^{n_1 k} - S^{n_1 k}$.  Also $p^{v_1 + 1} \parallel R^{n_1 p} - S^{n_1 p}$ (recall $p^{v_1}>2$).  By Lemma 3, $p \mid m$ for any $m$ such that $p^{v_1 + 1} \mid R^{n_1 m} - S^{n_1 m}$.  Now consider the binomial expansion of $((R^{n_1 p} - S^{n_1 p}) + S^{n_1 p})^{k}$ to see that $k=p$ is the least value of $k$ such that $p^{v_1+2} \mid R^{n_1 p k} - S^{n_1 p k}$.  Also $p^{v_1+2} \parallel R^{n_1 p^2} - S^{n_1 p^2}$.  By Lemma 3, $p^2 \mid m$ for any $m$ such that $p^{v_1+2} \mid R^{n_1 m} - S^{n_1 m}$.  Continuing in this way, we find that $p^d \mid m$ for any $m$ such that $p^{v_1 + d} \mid R^{n_1 m} - S^{n_1 m}$, giving the lemma.
\end{proof}

\begin{Lemma}  
Let $n_1$ and $n_2$ be positive integers with $n_1 \mid n_2$, and let $2^v \parallel \frac{n_2}{n_1}$.  Let $R$ and $S$ be be relatively prime odd positive integers with $R>S$, and let $2^t \parallel R^{n_1} - S^{n_1}$ and $2^u \parallel R^{n_1} + S^{n_1}$, with $h = \max(t,u)$.  Then, if $\frac{n_2}{n_1}$ is odd, $2^t \parallel R^{n_2} - S^{n_2}$ and $2^u \parallel R^{n_2} + S^{n_2}$.  And, if $\frac{n_2}{n_1}$ is even, $2^{h+v} \parallel R^{n_2} - S^{n_2}$ and $2 \parallel R^{n_2} + S^{n_2}$.  
\end{Lemma}

\begin{proof} 
If $\frac{n_2}{n_1}$ is odd, then $\frac{R^{n_2} - S^{n_2}}{R^{n_1} - S^{n_1}}$ and $\frac{R^{n_2} + S^{n_2}}{R^{n_1} + S^{n_1}}$ are both odd, so $2^t \parallel R^{n_2} - S^{n_2}$ and $2^u \parallel R^{n_2} + S^{n_2}$. So the lemma holds for $v=0$.   

Suppose $v=1$.  Then $n_2 = 2 m n_1$ for some odd $m$.  $2^{h+1} = 2^{t+u} \parallel (R^{m n_1} - S^{m n_1})(R^{m n_1} + S^{m n_1}) = R^{n_2} - S^{n_2}$ and (by consideration modulo 4) $2 \parallel R^{n_2} + S^{n_2}$.  So the lemma holds for $v=1$.  

The lemma follows by induction on $v$.  
\end{proof}

The following lemma adds a further restriction on one of the parameters in the infinite family (iv) in the Introduction.

\begin{Lemma}  
We cannot have $d = 1$ in (iv).
\end{Lemma}

\begin{proof} 
We consider the restrictions on the variables in (iv).  If $d=1$, then $h=2$, which requires $k-v = 1$, giving $k=2$. But then $(d+2)^k - d^k = 3^2 - 1 = 8$ has no odd divisor, contradicting $g>1$ and $w>0$.   
\end{proof}

\section{Preliminary Propositions}  

In this section and in the following section, we treat $a$, $b$, $c$ with $\gcd(a,b) > 1$, and define $Q$, $\alpha_p$, $\beta_p$, $\gamma_p$, $a_1$, $b_1$,  $c_1$ as in the first paragraph of Section 2.  Types A, B, C, and O are defined after Lemma 1.  

For any solution $(x,y,z)$ to (1.1) which is of Type A for a given prime $p \in Q$, we define $n = n(x,y,z)$ as follows:  letting 
$$ \frac{\beta_p}{ \gamma_p} = \frac{t}{s}, \gcd(s,t)=1, \eqno{(3.1)} $$
$n$ is the positive integer such that 
$$ y=ns, z=nt. $$
The existence of such $n$ follows from Lemma 1.  
For a specific solution $(x_i, y_i, z_i)$ we write $n(x_i,y_i,z_i) = n_i$, so that 
$$ y_i = n_i s,  z_i = n_i t.  \eqno{(3.2)} $$ 

We now define a function $f(n)$ relating to Type A solutions:   

For a given $(a,b,c)$, let $(x,y,z)$ be any solution to (1.1) which is Type A for some prime $p \in Q$, and let $s$ and $t$ be as in (3.1).    

Let $B$ be the set of all primes $q \in Q$ such that $\frac{\beta_q}{ \gamma_q} > \frac{t}{s}$.
Then the solution $(x, y, z)$ is Type B for every prime in $B$. 

Let $C$ be the set of all primes $q \in Q$ such that  $\frac{\beta_q}{ \gamma_q} < \frac{t}{s}$.  Then the solution $(x, y, z)$ is Type C for every prime in $C$.  

Let $A$ be the set of all primes $q \in Q$ for which $\frac{\beta_q}{ \gamma_q} = \frac{t}{s}$.  Then the solution $(x, y, z)$ is either Type A or Type O for every prime in $A$.  

Now we observe that any solution to (1.1) which is Type A for $p$ can be written as follows, noting that for some $n$, $y = n s$ and $z = nt$:  
$$ a_1^x \prod_{q \in Q} q^{\alpha_q x}  + b_1^{ns} \prod_{q \in Q} q^{\beta_q ns} = c_1^{nt} \prod_{q \in Q} q^{\gamma_q nt}.  \eqno{(3.3)} $$
The greatest number dividing all terms in (3.3) is 
$$ D = \prod_{q \in Q} (q^{\min(s \beta_q, t \gamma_q)})^n. \eqno{(3.4)} $$
Dividing both sides of (3.3) by $D$ and rearranging terms, we see that any solution of Type A for $p$ to (1.1) is equivalent to: 
$$(c_1^t \prod_{q \in C} q^{\gamma_q t - \beta_q s} )^n - (b_1^{s} \prod_{q \in B} q^{\beta_q s - \gamma_q t})^n = a_1^x \prod_{q \in A} q^{\alpha_q x - \gamma_q n t} .  \eqno{(3.5)}$$
Note that all variables in (3.5) other than $n$ itself are completely determined by $(a,b,c,p)$ except for $x$, which is determined by $n$.  We define a function $f(n)$: 
$$ f(n) = (c_1^t \prod_{q \in C} q^{\gamma_q t - \beta_q s} )^n - (b_1^{s} \prod_{q \in B} q^{\beta_q s - \gamma_q t})^n. $$
When $f(n)$ gives a solution to (1.1) as in (3.5), we have 
$$f(n) = a_1^x \prod_{q \in A} q^{\alpha_q x - \gamma_q n t}.$$  
For a specific solution $(x_i,y_i,z_i)$ we take $n_i$ as in (3.2) and write 
$$f(n_i) = a_1^{x_i} \prod_{q \in A} q^{\alpha_q x_i - \gamma_q n_i t}.$$ 

As in Section 2, we let $P(m)$ be the set of primes which divide a positive integer $m$.

\begin{Proposition} 
If (1.1) has two solutions $(x_1, y_1, z_1)$ and $(x_2,y_2,z_2)$ both of Type A for some prime $p \in Q$ with $P(f(n_1)) = P(f(n_2))$, then $N(a,b,c)=2$, and $(a,b,c, x_1, y_1, z_1, x_2, y_2, z_2)$ is in the infinite family (i) in the Introduction.  
\end{Proposition}

\begin{proof} 
Assume (1.1) has two distinct solutions $(x_1, y_1, z_1)$ and $(x_2,y_2,z_2)$ both of Type A for some prime $p \in Q$ for which (3.1) holds, let $n_1$, $n_2$ be as in (3.2) for $i \in \{1,2\}$, and further assume $P(f(n_1)) = P(f(n_2))$.  Since $(x_1,y_1,z_1)$ and $(x_2,y_2,z_2)$ are distinct, we have $n_1 \ne n_2$ by (3.2), so we can take 
$$n_1<n_2.$$  

Since $z_1 < z_2$, by Definition 2 in the Introduction the solutions $(x_1,y_1,z_1)$ and $(x_2,y_2,z_2)$ do not correspond to each other and are therefore not considered the same solution under Criterion 1, so that 
$$ N(a,b,c) > 1.$$

Write 
$$R = c_1^t \prod_{q \in C} q^{\gamma_q t - \beta_q s}, S = b_1^{s} \prod_{q \in B} q^{\beta_q s - \gamma_q t}.  $$  
Note that $\gcd(R,S)=1$. $R>S$ by (3.5).  We have
$$ R^n - S^n = f(n). $$  
Now we can apply Lemma 4 to $f(n_1)$ and $f(n_2)$ to obtain   
$$ R + S = 2^h, h \ge 2; n_1 = 1, n_2 = 2; f(n_1) = 2 j, 2 \nmid j; f(n_2) = 2^{h+1} j. \eqno{(3.6)}$$
Combining (3.6) with the right side of (3.5), we see that $a_1$ can have no odd divisors greater than one ($x_1 \ne x_2$ since $n_1 \ne n_2$ and $x$ increases with $n$ in (3.3)).  

Let $T = Q \cup \{ 2 \}$.  For every $q \in T$, let $W_q = \frac{w_2}{w_1}$ where $q^{w_1} \parallel a^{x_1}$ and $q^{w_2} \parallel a^{x_2}$.  We must have $W_q = \frac{x_2}{x_1} $ for every $q \in T$.  In considering $W_q$ we use $a^x = D f(n)$ for any $n$ which gives a solution to (1.1), where $D$ is as in (3.4).  For any odd $q \in A$, let $q^d \parallel j$, so that $W_q = \frac{n_2 t \gamma_q +d}{n_1 t \gamma_q + d} = \frac{2 t \gamma_q + d}{t \gamma_q + d} \le 2$.  For any odd $q \in B$, $W_q = \frac{n_2 t \gamma_q }{n_1 t \gamma_q } = 2$.  For any odd $q \in C$, $W_q = \frac{ n_2 s \beta_q}{n_1 s \beta_q} = 2$.  If $2 \in A$, then $W_2 = \frac{n_2 t \gamma_2 + h + 1}{n_1 t \gamma_2 + 1 } > 2$.  If $ 2 \mid a_1$, then $W_2 = h+1 > 2$.  Since we must have either $2 \in A$ or $2 \mid a_1$, we see, from the results on $W_q$, that $B$ and $C$ are empty and $A = \{2\}$ with $2 \nmid a_1$ since, by hypothesis, $A$ is not empty.  Recalling we have shown $a_1$ has no odd divisors, we have $a_1 = 1$, and, since $A = \{ 2\}$, from (3.6) we have $j = 1$.     

So we must have $p=2$, $a_1 = 1$, $Q = \{ 2 \}$, $R = c_1^t$, and $S = b_1^s$.   Since $n_1= j = 1$, from (3.6) we obtain $c_1^t - b_1^s = 2$ and $c_1^t + b_1^s = 2^h$, so that $b_1^s = 2^{h-1}-1$ and $c_1^t = 2^{h-1} + 1$.  If $h > 2$, it is a familiar elementary result that we must have $t=s=1$ except when $h=4$, $c_1 = 3$, $t =2$.  If $h=2$ then $b_1=1$ and we have the special case $(a_1, b_1, c_1) = (1,1,3)$, $Q = \{2\}$.  

So we must have one of the following:
$$ (a_1, b_1, c_1) = (1,2^{h-1}-1, 2^{h-1}+1), Q = \{2\}, h > 2, \eqno{(3.7)}$$
$$ (a_1, b_1, c_1) = ( 1, 7, 3), Q = \{2\},  \eqno{(3.8)}$$
$$ (a_1, b_1, c_1) = (1,1,3), Q = \{2\}.  \eqno{(3.9)}$$
These three cases are the only possibilities when (1.1) has two solutions $(x_1,y_1,z_1)$ and $(x_2,y_2,z_2)$ both of Type A for some prime in $Q$ with $P(f(n_1)) = P(f(n_2))$. We first show that for each of these three cases $N(a,b,c)=2$: since $N(a,b,c)>1$, it suffices to show that either no further solutions $(x,y,z)$ exist (which will handle (3.7) and (3.8)) or any further solutions $(x,y,z)$ must be considered the same solution as one of $(x_1,y_1,z_1)$ or $(x_2,y_2,z_2)$ under Criterion 1 (which will handle (3.9)).  

We first treat (3.7).  Writing $\alpha$ for $\alpha_2$, $ \beta$ for $\beta_2$, and $\gamma$ for $\gamma_2$, we have $\beta = \gamma$ since $s=t=1$.  So we can assume 
$$(a,b,c) = (2^\alpha, 2^\beta (2^{h-1} - 1), 2^\beta (2^{h-1} + 1)). \eqno{(3.10)}$$
We have the two solutions $(x_1,y_1,z_1)$ and $(x_2,y_2,z_2)$ with $n_1 = 1$ and $n_2=2$ as in (3.6), which must be the only solutions which are Type A for 2, since, if there existed a third such solution for some $n_3$, then, since $A = Q = \{ 2 \}$, we would have $P(R^{n_3} - S^{n_3}) = P(R^{n_2} - S^{n_2}) = P(R^{n_1} - S^{n_1})$, contradicting Lemma 4.  So it suffices to show that, for $(a,b,c)$ as in (3.10), (1.1) has no solutions of Type B, C, or O for 2. 

Since $a_1$, $b_1$, $c_1$ are all odd, there can be no solutions of Type O for 2.   We need to show there are no solutions of Type B or C for 2.  

For a solution of Type B for 2 we must have 
$$ 2^{z \beta} + 2^{y \beta} b_1^y = 2^{z \beta} (b_1 +2)^z  $$
so that 
$$ 2^{(y-z)\beta} b_1^y = (b_1 + 2)^z - 1, y>z, (y-z)\beta \ge h-1.  \eqno{(3.11)}$$
From (3.11) we have 
$$ (h-1) \log(2) + z \log(b_1) <  z \log(b_1+2) = z \log(b_1(1+\frac{2}{b_1})) < z \log(b_1) + \frac{2z}{b_1}. \eqno{(3.12)}$$
By Lemma 4, $z$ is the only value of $u$ such that $\rad((b_1+2)^u - 1) = \rad(2 b_1)$, so that, by Lemma 3, $z$ must be the least value of $u$ such that $b_1 \mid (b_1+2)^u - 1$,  so that $z \le \varphi(b_1) \le b_1 - 1$.  So now (3.12) gives 
$$ (h-1) \log(2) < \frac{2 (b_1-1)}{b_1} < 2 $$
which requires $h=3$ (since we have excluded $h=2$ from consideration in (3.7)).  

For the case $h=3$, (3.11) becomes
$$ 2^{(y-z)\beta} 3^y = 5^z - 1 $$
which, by Lemma 4, requires $z=2$, giving $y=1$, contradicting $y>z$.  So (3.7) does not allow a solution of Type B for 2.  

And a solution of Type C for 2 requires 
$$ 2^{y \beta} + 2^{y \beta} b_1^y = 2^{z \beta} (b_1+2)^z $$ 
so that 
$$ 1 + b_1^y = 2^{(z-y)\beta} (b_1+2)^z $$
which is impossible since $z>y$.  

Thus there is no third solution $(x_3,y_3,z_3)$ to (1.1) for the case (3.7), so $N(a,b,c)=2$ in this case.

Now we treat (3.8) for which a solution which is Type A for 2 requires, for some $v \in \intZ^+$, 
$$ 2^v + 7^y = 3^z. \eqno{(3.13)}$$
Consideration modulo 7 requires $2 \mid z$ in (3.13), so, if we assume there are no solutions of Type B, Type C, or Type O for 2, we can take $(a_1,b_1,c_1)= (1, 7, 9)$ and use (3.7), which we have shown gives only two solutions $(x,y,z)$.    

We need to show there are not solutions of Type B, Type C, or Type O for 2.  

A solution of Type O for 2 for the case (3.8) is impossible since $a_1=1$, $b_1=7$, and $c_1=3$ are all odd.

A solution which is Type B for 2 for the case (3.8) requires, for some $v \in \intZ^+$, 
$$ 1 + 2^v 7^y = 3^z. \eqno{(3.14)}$$
But then we note that $7 \mid 3^z-1$ requires $13 \mid 3^z-1$, a contradiction which shows that there can be no solution of Type B for 2.   

A solution of Type C for 2 for the case (3.8) requires, for some $v \in \intZ^+$, 
$$ 1 + 7^y = 2^v 3^z \eqno{(3.15)}$$
which is impossible modulo 3.

Thus there is no third solution $(x_3,y_3,z_3)$ to (1.1) for the case (3.8), so $N(a,b,c)=2$ in this case.  

Now we treat the case (3.9), letting $\alpha = \alpha_2$, $\beta = \beta_2$, and $\gamma=\gamma_2$.  In this case we have, for any solution $(x,y,z)$,  
$$ 2^{\alpha x} + 2^{\beta y} = 2^{\gamma z} 3^z, $$
making solutions of Type C for 2 and Type O for 2 clearly impossible.  

For a solution which is Type A for 2, we have  
$$ 2^{\alpha x} + 2^{\gamma z} = 2^{\gamma z} 3^z \eqno{(3.16)} $$
so that for a given value of $z$ there is at most one solution of Type A for 2. 

Similarly, a solution of Type B for 2 requires 
$$ 2^{\gamma z} + 2^{\beta y} = 2^{\gamma z} 3^z, \eqno{(3.17)} $$ 
so that for a given value of $z$ there is at most one solution of Type B for 2.    

Suppose that for a given choice of $z$ we have both a solution which is Type A for 2 and a solution which is Type B for 2: letting $(x_a, y_a, z)$ be the solution which is Type A for 2 in (3.16) and letting $(x_b, y_b, z)$ be the solution which is Type B for 2 in (3.17), we have $\alpha x_a = \beta y_b$, and we see that $a^{x_a} = b^{y_b}$ and $a^{x_b} = b^{y_a}$.  Using Criterion 1 we see that $(x_a, y_a, z)$ and $(x_b,y_b,z)$ are to be considered the same solution.  So, to obtain $N(a,b,c)=2$, it suffices to show that there are only two values of $z$ which are possible in solutions to (1.1) which are Type A for 2.  

From (3.16) we see that any solution which is Type A for 2 must give 
$$ 3^z - 1 = 2^{\alpha x - \gamma z}, $$
which requires either $z=1$ or $z=2$ by Lemma 4.  Thus (using Criterion 1) we see that $N(a,b,c)=2$ for the case (3.9).  

It remains to show that $(a,b,c, x_1, y_1, z_1, x_2, y_2, z_2)$ is in the infinite family (i) in the Introduction.  We have shown that if (1.1) has two solutions $(x_1,y_1,z_1)$ and $(x_2,y_2,z_2)$ both of Type A for some prime in $Q$ with $P(f(n_1))=P(f(n_2))$, then one of (3.7), (3.8), or (3.9) must hold.  We have also shown that $(x_1,y_1,z_1)$ and $(x_2,y_2,z_2)$ must satisfy (3.6), with $j=a_1=1$.  If (3.7) holds, $s=t=1$, so, using (3.2) and the right side of (3.5) in combination with (3.6) and the definition of $f(n)$, we find that when (3.7) holds we must have $(a,b,c)$ as in (3.10) with $(x_1,y_1,z_1)$ and $(x_1,y_2,z_2)$ as follows: 
$$  (x_1, y_1, z_1)= (\frac{\beta+1}{\alpha}, 1, 1), (x_2,y_2,z_2) = (\frac{2 \beta + h+1}{\alpha}, 2, 2).  $$
These two solutions correspond to two solutions to (1.1) for $(a,b,c) = (A, B, C)$ as follows:
$$ (A,B,C) = (2, 2^{\beta} (2^{h-1}-1), 2^\beta (2^{h-1}+1)), (X_1, Y_1, Z_1) = (\beta+1, 1, 1), (X_2, Y_2, Z_2) = (2 \beta + h + 1, 2,2).  \eqno{(3.18)}$$
We see that $(A,B,C, X_1, Y_1, Z_1, X_2,Y_2,Z_2)$ is a member of $F$ for the infinite family (i).  Also we have shown that, when (3.8) holds, $(x_1,y_1,z_1)$ and $(x_2,y_2,z_2)$ must  correspond to solutions in which $(1,7,3)$ is replaced by $(1,7,9)$, so again we have (3.7) giving (3.18).  Finally, if (3.9) holds, we have shown that $(x_1,y_1,z_1)$ and $(x_2,y_2,z_2)$ must satisfy $z_1=1$ and $z_2=2$, so that $\beta_2 y_1 = \gamma_2$ and $\beta_2 y_2 = 2 \gamma_2$.  Letting $u= \beta_2 y_1$, we have a member of $F$ (for the infinite family (i)) whose two solutions correspond to $(x_1,y_1,z_1)$ and $(x_2,y_2,z_2)$.  By Definition 3, in all three cases $(a,b,c, x_1,y_1,z_1,x_2,y_2,z_2)$ is in the infinite family (i).   
\end{proof}

\begin{Corollary3_1*}  
If (1.1) has two solutions $(x_1,y_1,z_1)$ and $(x_2,y_2,z_2)$ both of which are Type A for some prime in $Q$ and if further (1.1) has no solutions of Type O for any prime in $Q$, then $N(a,b,c)=2$, and $(a,b,c, x_1, y_1, z_1, x_2, y_2, z_2)$ is in the infinite family (i) in the Introduction.  

This statement holds if we replace \lq\lq Type A'' by \lq\lq Type B''. 
\end{Corollary3_1*}

\begin{proof} 
By the symmetry of $a$ and $b$ (using Definitions 2 and 3) it suffices to prove the first paragraph of the Corollary.  

Choose a prime $p \in Q$ such that $(x_1,y_1,z_1)$ and $(x_2,y_2,z_2)$ are both Type A for $p$, and let $s$ and $t$ be as in (3.1) and (3.2) with $i \in \{1,2\}$.  We have (3.5) with $(n,z) = (n_i, z_i)$, $i \in \{1,2\}$, where $n_1 \ne n_2$ and the sets $A$, $B$, and $C$ are derived from (3.1).  Since for the case under consideration no solutions of Type O are possible for any prime in $Q$, we see that $(x_1,y_1,z_1)$ and $(x_2,y_2,z_2)$ are both Type A for every prime in $A$, so, for every prime $q \in A$, we have 
$$ \alpha_q x_1 - \gamma_q n_1 t > 0, \alpha_q x_2 - \gamma_q n_2 t > 0$$
so that, using (3.5) and recalling the definition of $f(n)$, we have $P(f(n_1)) = P(f(n_2))$, so we can apply Proposition 3.1 to obtain the corollary.
\end{proof}

\begin{Proposition}  
Assume (1.1) has a solution which is Type C for some prime $p \in Q$.  Then (1.1) can have no solution which is Type A, B, or O for $p$.   
\end{Proposition}

\begin{proof}
Let $p$ be any prime in $Q$.  For brevity, write $\alpha$ for $\alpha_p$, $\beta$ for $\beta_p$, $\gamma$ for $\gamma_p$. Let $(x_c, y_c, z_c)$ be a solution to (1.1) which is Type C for $p$ and let $(x_d,y_d,z_d)$ be a solution to (1.1) which is Type A, B, or O for $p$.  By Lemma 1 $\alpha x_c < \gamma z_c$ and $\alpha x_d \ge \gamma z_d$ so that
$$ \frac{x_c}{z_c} < \frac{x_d}{z_d}. \eqno{(3.19)}$$
Similarly
$$ \frac{y_c}{z_c} < \frac{y_d}{z_d}. \eqno{(3.20)}$$
Now let $a_p = a/p^\alpha$, $b_p = {b}/{p^\beta}$ and let $c_p = {c}/{p^\gamma}$.  $a^{x_d} < c^{z_d}$ so $p^{\alpha x_d} a_p^{x_d} < p^{\gamma z_d} c_p^{z_d}$, so that $ p^{\alpha x_d - \gamma z_d} a_p^{x_d} < c_p^{z_d} $, so that
$$ x_d \log(a_p) < z_d \log(c_p).  \eqno{(3.21)}$$
Similarly, 
$$ y_d \log(b_p) < z_d \log(c_p).  \eqno{(3.22)}$$
Now suppose $a^{x_c} \ge b^{y_c}$. Then $2 a^{x_c} \ge c^{z_c}$, so that $2 a_p^{x_c} \ge p^{\gamma z_c - \alpha x_c} c_p^{z_c} \ge 2 c_p^{z_c}$, so that 
$$ x_c \log(a_p) \ge z_c \log(c_p). \eqno{(3.23)}$$
Since (3.21) requires $c_p > 1$, (3.23) gives $a_p > 1$.  So (3.23) in combination with (3.21) gives 
$$\frac{x_c}{z_c} > \frac{x_d}{z_d}, \eqno{(3.24)}$$
contradicting (3.19).   So we must have $b^{y_c} > a^{x_c}$, so that $2 b^{y_c} > c^{z_c}$, and, recalling (3.22), the same argument which handles the case $a^{x_c} > b^{y_c}$ yields 
$$ \frac{y_c}{z_c} > \frac{y_d}{z_d}, \eqno{(3.25)}$$
contradicting (3.20).  So the existence of a solution of Type C for $p$ makes solutions of Type A, B, or O for $p$ impossible.  
\end{proof}

\begin{Proposition}  
For a given $p \in Q$, there is at most one solution $(x,y,z)$ to (1.1) which is Type O for $p$.  
\end{Proposition}

\begin{proof} 
Let $r$, $s$, and $t$ be positive integers with $\gcd(r,s,t)=1$ such that
$$ r \alpha_p = s \beta_p = t \gamma_p, p \in Q. $$
By Lemma 1 any solution $(x,y,z)$ which is Type O for $p$ must have 
$$ x = n r, y = n s, z = n t, $$
for some positive integer $n$.  
By Fermat's Last Theorem \cite{W}, $n \le 2$. So if there are two distinct solutions $(x,y,z)$ to (1.1) which are Type O for $p$, then $a^r + b^s = c^t$ and $a^{2 r} + b^{2s} = c^{2t}$, which is impossible since $c^{2t} = (a^r + b^s)^2 > a^{2r} + b^{2s} = c^{2t}$.  
\end{proof}

Let $S_c$ be the set of all triples $(a,b,c)$ such that $\gcd(a,b)>1$, $N(a,b,c)>1$, $a$, $b$, and $c$ are not all powers of 2, and there exists a solution to (1.1) which is Type C for some prime in $Q$.

\begin{Proposition}  
For $(a,b,c) \in S_c$, (1.1) has exactly two solutions $(x_1, y_1,z_1)$ and $(x_2,y_2,z_2)$, $N(a,b,c)=2$, and $(a,b,c,\allowbreak  x_1, y_1,z_1,\allowbreak  x_2, y_2, z_2)$ is in the infinite family (ii).    
\end{Proposition}

\begin{proof} 
After Proposition 3.2, it suffices to show that if (1.1) has two solutions $(x_1,y_1,z_1)$ and $(x_2, y_2,z_2)$ both of which are Type C for some prime in $Q$, then these two solutions are the only solutions $(x,y,z)$ to (1.1), $N(a,b,c)=2$, and $(a,b,c,x_1,y_1,z_1, x_2, y_2,z_2)$ is in the infinite family (ii).  

We define a function $f(n)$ relating to Type C solutions:  

For a given $(a,b,c)$, let $(x,y,z)$ be any solution to (1.1) which is Type C for some prime $p \in Q$, and let $r$ and $s$ satisfy
$$ \frac{\alpha_p }{\beta_p} = \frac{s}{r}, \gcd(r,s)=1.  \eqno{(3.26)} $$
Then by Lemma 1 we must have a positive integer $n = n(x,y,z)$ such that  
$$ x = n r, y = n s.  \eqno{(3.27)}$$

Let $A$ be the set of all primes $q \in Q$ such that $\frac{\alpha_q}{\beta_q} > \frac{s}{r}$.  Then the solution $(x, y, z)$ is Type A for every prime in $A$.  

Let $B$ be the set of all primes $q \in Q$ such that $\frac{\alpha_q}{\beta_q} < \frac{s}{r}$.  Then the solution $(x,y, z)$ is Type B for every prime in $B$.   

Let $C$ be the set of all primes $q \in Q$ for which $\frac{\alpha_q}{\beta_q} = \frac{s}{r}$.  Then the solution $(x,y,z)$ is either Type C or Type O for every prime in $C$.   

Now we observe that any solution to (1.1) which is Type C for $p$ can be written as follows, using (3.27):  
$$ a_1^{nr} \prod_{q \in Q} q^{\alpha_q n r} + b_1^{ns} \prod_{q \in Q} q^{\beta_q n s} = c_1^{z} \prod_{q \in Q} q^{\gamma_q z}. \eqno{(3.28)}$$
The greatest number dividing all terms in (3.28) is 
$$ D  = \prod_{q \in Q} (q^{\min(r \alpha_q, s \beta_q)})^n.  \eqno{(3.29)}$$
Dividing both sides of (3.28) by $D$, we see that any solution of Type C for $p$ to (1.1) is equivalent to:  
$$(a_1^r \prod_{q \in A} q^{\alpha_q r - \beta_q s})^n + ( b_1^s \prod_{q \in B} q^{\beta_q s - \alpha_q r} )^n = c_1^z \prod_{q \in C} q^{\gamma_q z - \alpha_q n r}. \eqno{(3.30)}$$
Note that all variables in (3.30) other than $n$ itself are completely determined by $(a,b,c,p)$ except for $z$, which is determined by $n$.  
We define a function $f(n)$:
$$ f(n) = (a_1^r \prod_{q \in A} q^{\alpha_q r - \beta_q s})^n + ( b_1^s \prod_{q \in B} q^{\beta_q s - \alpha_q r} )^n.$$
When $f(n)$ gives a solution to (1.1), 
$$f(n) = c_1^z \prod_{q \in C} q^{\gamma_q z - \alpha_q n r}.$$

Let 
$$R = a_1^r \prod_{q \in A} q^{\alpha_q r - \beta_q s}, S = b_1^s \prod_{q \in B} q^{\beta_q s - \alpha_q r}. \eqno{(3.31)}$$
Note $\gcd(R,S)= 1$.  Assume first we do not have the special case $R=S=1$.  Then it suffices to prove Proposition 3.4 for the case $R>S$, by the symmetry of $a$ and $b$ (recall Definitions 2 and 3). So we have 
$$f(n) = R^n + S^n, R>S, \gcd(R,S)=1. $$  

Now fix $(a,b,c)$ and assume (1.1) has two distinct solutions $(x_1, y_1, z_1)$ and $(x_2, y_2, z_2)$ both of which are Type C for a given $p \in Q$ with $r$ and $s$ as in (3.26), so that by (3.27) 
$$ x_i = n_i r, y_i = n_i s, i \in \{1,2\}, n_i \in \intZ^+, $$
so that (3.30) holds for $(n,z) = (n_1, z_1)$ and for $(n,z) = (n_2, z_2)$.  
Since $(x_1,y_1,z_1)$ and $(x_2,y_2,z_2)$ are distinct, we can take 
$$n_1 < n_2,$$
so, by Definition 2, clearly $(x_1,y_1,z_1)$ and $(x_2,y_2,z_2)$ do not correspond to each other.  

By Propositions 3.2 and 3.3, both $(x_1,y_1,z_1)$ and $(x_2,y_2,z_2)$ must be Type C for every prime $q \in C$, so, for every such $q$,  we have $\gamma_q z - \alpha_q  n r > 0$ in (3.30) for $(n,z)=(n_i,z_i)$, $i \in \{ 1,2\}$, so that   
$P(f(n_1)) = P(f(n_2))$ where, as before, $P(m)$ is the set of primes which divide $m$.   
Now we can apply Lemma 5 to $f(n_1)$ and $f(n_2)$ to obtain   
$$ R = 2, S=1, n_1 =1, n_2 = 3; f(n_1) = 3, f(n_2) = 9. \eqno{(3.32)} $$
Noting that (3.32) is obtained under no assumptions other than the assumption that  $(x_1,y_1,z_1)$ and $(x_2,y_2,z_2)$ are both Type C for the given prime $p$, we see that $(x_1, y_1,z_1)$ and $(x_2,y_2,z_2)$ must be the only solutions to (1.1) of Type C for this $p$.  Since by Proposition 3.2 we need only consider solutions of Type C for this $p$, we see that (1.1) has exactly the two solutions $(x_1,y_1,z_1)$ and $(x_2,y_2,z_2)$ when $R > S$.  Since these two solutions do not correspond to each other, we have $N(a,b,c)=2$.  To complete the treatment of the case $R>S$, it remains to show that $(a,b,c,x_1,y_1,z_1, x_2,y_2,z_2)$ is in the infinite family (ii).  
 
Since we are assuming there is at least one prime $p \in Q$ for which the solutions $(x_1, y_1, z_1)$ and $(x_2,y_2, z_2)$ are Type C for $p$, from (3.32) we see that $C = \{3\}$ and $c_1 = 1$.  Also from (3.32) we have $b_1=1$ with $B$ empty.  

Assume $A = \{ 2 \}$. Then, by (3.31) and(3.32), $a_1 = b_1 = c_1 = 1$, $Q = \{ 2,3\}$, and, since $R+S = 2 + 1 = 3$ and $R^3 + S^3 = 2^3 + 1 = 3^2$, both solutions are Type A for 2 and neither solution is Type O for either prime in $Q$. By the Corollary to Proposition 3.1, we must have $(a,b,c, x_1, y_1, z_1, x_2, y_2, z_2)$ in the infinite family (i) in the Introduction.  But then $Q = \{ 2 \}$, contradicting $Q = \{2,3\}$.   (Note that $Q = \{2\}$ would require $C$ empty and $c_1 = 3$, giving (3.9)). So $A \ne \{2\}$.  

So, by (3.31) and (3.32), we must have $a_1^r =2$ with $A$ empty and $Q = \{ 3\}$. 
Thus we have $a_1 = 2$, $b_1 = 1$, $c_1 = 1$. So we have 
$$ a = 2 \cdot 3^{\alpha_3}, b = 3^{\beta_3}, c = 3^{\gamma_3} $$
with exactly the two solutions $(x_1,y_1,z_1)$ and $(x_2,y_2,z_2)$.  Since $a_1^r = 2$, we have $r=1$, so (3.27) and (3.32) give $x_1=1$, $x_2=3$.  Since $(x_1,y_1,z_1)$ is Type C for 3, we have $\beta_3 y_1 = \alpha_3$, and, using (3.32) and noting $f(n_1) = 3^{\gamma_3 z_1 - \alpha_3}$, we have $\gamma_3 z_1 = \alpha_3 + 1$.  So $(x_1,y_1,z_1) = (1, \alpha_3/\beta_3, (\alpha_3+1)/\gamma_3)$.  Similarly, $(x_2,y_2,z_2) = (3, 3 \alpha_3 / \beta_3, (3 \alpha_3 + 2)/\gamma_3)$.  These two solutions correspond to the solutions $(1,\alpha_3, \alpha_3+1)$ and $(3, 3 \alpha_3, 3 \alpha_3 + 2)$ to (1.1) for $(a,b,c) = (A,B,C)$ as follows:     
$$ A = 2 \cdot 3^{\alpha_3}, B = 3, C = 3.  \eqno{(3.33)} $$
Since $(A,B,C)$ as in (3.33) and its two solutions give a nine-tuple in $F$ for the infinite family (ii), we see that, by Definition 3, Proposition 3.4 holds for the case $R>S$.      

If $R=S=1$, then, by (3.31), $a_1 = b_1 = 1$ and $A$ and $B$ are empty, so that $f(n) = 2$ for every $n$.  Since we are assuming $C$ has at least one element, we have $c_1=1$ and $C = Q = \{2\}$.  So in this case all of $a$, $b$, and $c$ must be powers of 2, which has been excluded.  This completes the proof of Proposition 3.4.  
\end{proof}

\begin{Proposition}   
Suppose for some $(a,b,c)$ we have a nonempty subset $Q_1$ of $Q$ such that, for any primes $q_1$ and $q_2$ in $Q_1$, 
$$ \frac{\alpha_{q_1} }{\alpha_{q_2}} = \frac{\beta_{q_1}}{\beta_{q_2}} = \frac{\gamma_{q_1}}{\gamma_{q_2}}. \eqno{(3.34)} $$
Let $Q_2$  be the set of all primes in $Q$ which are not in $Q_1$.  
Then we can take
$$ a = a_1 g^{\alpha_g} \prod_{q \in Q_2} q^{\alpha_q}, b = b_1 g^{\beta_g} \prod_{q \in Q_2} q^{\beta_q}, c=  c_1 g^{\gamma_g} \prod_{q \in Q_2} q^{\gamma_q} \eqno{(3.35)} $$ 
where $\alpha_g$, $\beta_g$, and $\gamma_g$ are positive integers and $g$ is a positive integer divisible by every prime in $Q_1$ and by no other prime.  
\end{Proposition}

\begin{proof} 
We use the notation indicating proportions in which (3.34) would be written as follows:
$$\alpha_{q_1} : \alpha_{q_2} = \beta_{q_1} : \beta_{q_2} = \gamma_{q_1} : \gamma_{q_2}.$$
Assume $Q_1$ is composed of $n$ primes.  Then we can find a set of $n$ positive integers $t_1$, $t_2$, \dots, $t_n$ with no common divisor such that 
$$ \alpha_{q_1}: \alpha_{q_2} : \dots : \alpha_{q_n} =\beta_{q_1} : \beta_{q_2} : \dots : \beta_{q_n} = \gamma_{q_1} : \gamma_{q_2} : \dots : \gamma_{q_n} = t_1: t_2: \dots : t_n $$
and let 
$$ \alpha_{q_i} = h t_i, \beta_{q_i} = j t_i, \gamma_{q_i} = m t_i, 1 \le i \le n, \eqno{(3.36)}$$
for some positive integers $h$, $j$, and $m$.  

So we can take $g = q_1^{t_1} q_2^{t_2} \cdots q_n^{t_n}$, $\alpha_g = h$, $\beta_g = j$, and $\gamma_g = m$.    
\end{proof}

Since any solution to (1.1) of a given Type for some $q_1 \in Q_1$ is of the same Type for any $q_2 \in Q_1$, Proposition 3.5 allows us to refer to solutions of Type A for $g$, Type B for $g$, etc., even when $g$ is composite.  We say that a solution is of a given Type for $g$ when this solution is of that same Type for every prime dividing $g$, where the set of $n$ primes dividing $g$ satisfies (3.36).

\section{Proof of Theorem 1.1 and further results}   

Let $S_o$ be the set of all triples $(a,b,c)$ such that $\gcd(a,b)>1$, $N(a,b,c)>1$, and there exists a solution to (1.1) which is Type O for some prime in $Q$.

\begin{Proposition}  
For $(a,b,c) \in S_o$ we have $N(a,b,c)=2$, and (1.1) has two solutions $(x_1,y_1,z_1)$ and $(x_2,y_2,z_2)$ with $(a,b,c, \allowbreak  x_1,y_1,z_1,\allowbreak  x_2,y_2,z_2)$ in one of the infinite families (iii) or (iv).    
\end{Proposition}

\begin{proof}
Choose $(a,b,c) \in S_o$ and let $(x_1, y_1,z_1)$ be a solution to (1.1) which is Type O for $p$ where $p \in Q$.  Since $N(a,b,c)>1$, we can assume (1.1) has a second solution $(x_2,y_2,z_2)$.  By Propositions 3.2 and 3.3, $(x_2,y_2,z_2)$ must be either Type A or Type B for $p$.  By the symmetry of $a$ and $b$ (using Definitions 2 and 3) we can assume $(x_2,y_2,z_2)$ is Type A for $p$.  By the Corollary to Lemma 2, $(x_2,y_2,z_2)$ must be Type A for every prime in $Q$. 

Let 
$$\frac{\beta_p}{\gamma_p}   = \frac{t}{s}, \gcd(s,t)=1.$$
We can take      
$$ y_1 = n_1 s, z_1 = n_1 t, y_2 = n_2 s, z_2 = n_2 t, n_1 \in \intZ^+, n_2 \in \intZ^+, n_1 \ne n_2.$$
Since, for every prime $q \in Q$, $y_2 \beta_q = z_2 \gamma_q$, we see that $\frac{\beta_q}{\gamma_q} = \frac{t}{s}$ for every prime $q \in Q$, and thus $(x_1,y_1,z_1)$ must be either Type O or Type A for every prime in $Q$.

Let $G$ be the set of all primes in $Q$ for which $(x_1,y_1,z_1)$ is Type O and let $H$ be the set of all primes in $Q$ for which $(x_1, y_1, z_1)$ is Type A.  For every prime $q \in G$, $\alpha_q x_1 = \beta_q y_1 = \gamma_q z_1$, so that for any two primes $q_1$ and $q_2$ in $G$ we have 
$$   \frac{\alpha_{q_1}}{\alpha_{q_2}} = \frac{\beta_{q_1}}{\beta_{q_2}} = \frac{\gamma_{q_1}}{\gamma_{q_2}}$$
so that we can apply Proposition 3.5 to see that we can write 
$$ a = a_1 g^{\alpha_g} \prod_{q \in H} q^{\alpha_q}, b = b_1 g^{\beta_g} \prod_{q \in H} q^{\beta_q}, c=  c_1 g^{\gamma_g} \prod_{q \in H} q^{\gamma_q}  \eqno{(4.1)}$$
for some positive integers $g$, $\alpha_g$, $\beta_g$, $\gamma_g$, with $P(g) = G$.  Note $g > 1$ since $p \in G$.    

From the solution $(x_1,y_1,z_1)$ we derive 
$$ c_1^{n_1 t} - b_1^{n_1 s} = a_1^{x_1} \prod_{q \in H} q^{\alpha_q x_1 - n_1 t \gamma_q}, x_1 = \frac{n_1 t \gamma_g}{\alpha_g}  \eqno{(4.2)}$$
From the solution $(x_2,y_2,z_2)$ we derive  
$$ c_1^{n_2 t} - b_1^{n_2 s} = a_1^{x_2} g^{\alpha_g x_2 - n_2 t \gamma_g} \prod_{q \in H} q^{\alpha_q x_2 - n_2 t \gamma_q}, x_2 > \frac{n_2 t \gamma_g}{\alpha_g}. \eqno{(4.3)}  $$
From the expressions for $x_1$ and $x_2$ in (4.2) and (4.3) we derive 
$$ \frac{x_2}{x_1} > \frac{n_2}{n_1}. \eqno{(4.4)}$$

Let $R = c_1^{t}$, $S = b_1^{s}$, so that $R>S$.  Let 
$$ f(n) = c_1^{t n} - b_1^{s n} = R^n - S^n.  $$

Let $U$ be the product of all primes in $H$, let $r$ be any prime dividing $a_1 U$ (if such $r$ exists), and let $r^{v_1} \parallel f(n_1)$ and $r^{v_2} \parallel f(n_2)$.  Recalling (4.2) and (4.3) we see that, if $r \mid a_1$, then $\frac{v_2}{v_1} = \frac{x_2}{x_1} > \frac{n_2}{n_1}$ by (4.4).  If $r \in H$ then by (4.2), (4.3), and (4.4)   
$$\frac{v_2}{v_1} = \frac{\alpha_r x_2 - n_2 t \gamma_r}{\alpha_r x_1 - n_1 t \gamma_r} = \frac{\alpha_r \frac{x_2}{x_1} x_1 - \frac{n_2}{n_1} n_1 t \gamma_r} {\alpha_r x_1 - n_1 t \gamma_r}> \frac{n_2}{n_1}. $$ 
So in either case, we have 
$$ \frac{v_2}{v_1} > \frac{n_2}{n_1}. \eqno{(4.5)} $$

Let $n_0$ be the least value of $n$ such that $\rad(a_1 U) \mid f(n)$. By Lemma 3
$$ n_0 \mid n_1, n_0 \mid n_2  \eqno{(4.6)}$$
so that $f(n_0) \mid f(n_1)$, so that $f(n_0) $ can be divisible by no primes which do not divide $a_1 U$, so that
$$P(f(n_0)) = P(f(n_1)) = P(a_1 U).$$
Assume there do not exist distinct positive integers $m_1$ and $m_2$ such that $P(f(m_1)) = P(f(m_2))$.   Then $n_1$ must be the only value of $n$ such that $P(f(n)) = P(a_1 U)$, so that $n_0 = n_1$, and, by (4.6), $n_1 \mid n_2$.   

Now assume we have $m_1$ and $m_2$ such that $P(R^{m_1} - S^{m_1}) = P(R^{m_2} - S^{m_2})$.  Let $h$ be an integer such that $2^h \parallel R+S$.  By Lemma 4, $\{ m_1, m_2 \} = \{ 1, 2\}$, $R-S \equiv 2 \bmod 4$, and $R + S = 2^h$ with $h \ge 2$.  If $n_1 \not\in \{1, 2 \}$, then the above argument still applies to show $n_1 \mid n_2$.  If $n_1 = 1$, then again we have $n_1 \mid n_2$.  

So we consider $n_1 = 2$ with $P(f(1))=P(f(2))$: 
Let $2^{v_1} \parallel f(n_1) = f(2) = (R-S)(R+S)$ so that $v_1 = h + 1$.  By (4.2), $\gcd(f(n_1), g) = 1$, so, since $P(f(1)) = P(f(2))$ and $g \mid f(n_2)$, we must have $n_2 > 2 = n_1$, so that, letting $2^{v_2} \parallel f(n_2)$, by (4.5) we have $v_2 > v_1 = h+1 > 2$.  Since $n_1 = 2$ is the least value of $n$ such that $2^2 \mid f(n)$, by Lemma 3 we have $n_1 \mid n_2$.  

So we can assume 
$$ n_1 \mid n_2.  \eqno{(4.7)} $$

Now we can apply Lemma 6: letting $r$ be any prime dividing $a_1 U$, and letting $r^{v_1} \parallel f(n_1)$ and $r^{v_2} \parallel f(n_2)$, we see that, if $r^{v_1} > 2$, 
$$ r^{v_2-v_1} \mid \frac{n_2}{n_1} $$
so that by (4.5) 
$$\frac{n_2}{n_1} \ge r^{v_2-v_1} \ge r^{\frac{v_2}{v_1} - 1} > r^{\frac{n_2}{n_1} - 1} \eqno{(4.8)} $$
which is impossible.  So we must have $r^{v_1} \le 2$.  Since $r$ can be any prime dividing $a_1 U$ we see that $H$ contains no odd primes and $a_1 \le 2$, giving only three possible cases:

Case 1: $a_1 = 1$, $H = \{2\}$, $f(n_1) = 2$.

Case 2:  $a_1 = 2$, $H$ is empty, $f(n_1) = 2$.

Case 3:  $a_1 = 1$, $H$ is empty, $f(n_1) = 1$.   

For all three cases, $n_1 = 1$.  

For Case 1, we have, for $a$, $b$, $c$ as in (4.1),  
$$ a = 2^{\alpha_2} g^{\alpha_g}, b = 2^{\beta_2} g^{\beta_g} b_1, c = 2^{\gamma_2} g^{\gamma_g} c_1. \eqno{(4.9)} $$
Let 
$$ d = b_1^{y_1}, b_0 = b^{y_1},  c_0 = c^{z_1}.$$
Since $c_1^{z_1} - b_1^{y_1} = c_1^{n_1 t} - b_1^{n_1 s} = c_1^{t} - b_1^{s} = f(n_1) = 2$ we have
$$ c_1^{z_1} = d + 2.  \eqno{(4.10)} $$
So we have 
$$ b_0 = 2^{\beta_2 y_1} g^{\beta_g y_1} d, c_0 = 2^{\gamma_2 z_1} g^{\gamma_g z_1} (d + 2).  \eqno{(4.11)}$$
Since the solution $(x_1,y_1, z_1)$ is Type O for $g$, we have $\beta_g y_1 =  \gamma_g z_1  = \alpha_g x_1$; also since $(x_1,y_1,z_1)$ is Type A for 2 and $f(n_1) = 2$, we have $\beta_2 y_1 = \gamma_2 z_1 = \alpha_2 x_1 - 1$ by (4.2).  So from (4.11) we have  
$$ b_0 = 2^{\alpha_2 x_1 - 1} g^{\alpha_g x_1} d, c_0 = 2^{\alpha_2 x_1 - 1} g^{\alpha_g x_1} (d+2).  \eqno{(4.12)}$$
Let $k = n_2$, and note that, since $n_1 = 1$,  
$$ k = \frac{n_2 s}{s} = \frac{y_2}{y_1}, k = \frac{n_2 t}{t} = \frac{z_2}{z_1}.  $$
The solutions $(x_1,y_1,z_1)$ and $(x_2,y_2,z_2)$ to (1.1) for the triple $(a,b,c)$ as in (4.9) correspond to the solutions $(x_1, 1, 1)$ and $(x_2, k, k)$ for the triple $(a, b_0, c_0)$.  

Since $g^{\alpha_g x_2} \parallel a^{x_2} = c_0^{k} - b_0^{k}$, from (4.12) we find that $\alpha_g x_2 = k \alpha_g x_1 + w_g$ where $g^{w_g} \parallel (d + 2)^k - d^k$ (such $w_g>0$ must exist since $(x_2,y_2,z_2)$ is Type A for $g$; note that $g^{w_g}$ is the greatest odd divisor of $(d+2)^k-d^k$ and $\alpha_g$ must divide $w_g$).  So 
$$ x_2 = k x_1 + \frac{w_g}{\alpha_g}.  \eqno{(4.13)} $$

Since $2^{\alpha_2 x_2} \parallel a^{x_2} = c_0^{k} - b_0^{k}$, from (4.12) we find that $\alpha_2 x_2 = k \alpha_2 x_1 - k + w_2$ where $2^{w_2} \parallel (d+2)^k - d^k$.  
But also by (4.13) $\alpha_2 x_2 = k \alpha_2 x_1 + \frac{\alpha_2 w_g}{\alpha_g}$, so that $w_2 = k + \frac{\alpha_2 w_g}{\alpha_g} > 1$, so that consideration modulo 4 gives 
$$ 2 \mid k.$$  
So now we can use Lemma 7 to see that we must have $k + \frac{\alpha_2 w_g}{\alpha_g} = w_2 = h+v$ where $2^v \parallel k$ and $2^h \parallel 2d + 2$.  So we have 
$$k-v = h - \frac{\alpha_2 w_g}{\alpha_g}.$$
Now from (4.12) and (4.13) we see that $(a,b_0,c_0,x_1, 1, 1,\allowbreak k x_1 + ({w_g}/{\alpha_g}), k,k)$ (and therefore also $(a,b,c,\allowbreak x_1,y_1,z_1, \allowbreak x_2,y_2,z_2)$) is in the infinite family (iv) in the Introduction when $i = \alpha_2$, $j = \alpha_g$, $u = x_1$, $w = w_g$, and $iu > 1$.  To complete the treatment of Case 1, we need to show $N(a,b,c)=2$.  Since $(x_1,y_1,z_1)$ and $(x_2,y_2,z_2)$  do not correspond to each other, it suffices to show that there is no third solution $(x_3,y_3,z_3)$ for the triple $(a,b,c)$.  

First note that $d = b_1^{y_1}$ in the member of the infinite family (iv) derived from the solutions $(x_1,y_1,z_1)$ and $(x_2,y_2,z_2)$. By Lemma 8, $b_1 \ne 1$.  

Now suppose there is a third solution other than $(x_1,y_1,z_1)$ and $(x_2,y_2, z_2)$ for the triple $(a,b,c)$.  By Propositions 3.2 and 3.3 this third solution $(x_3,y_3,z_3)$ must be either Type A for $g$ or Type B for $g$.  If this third solution is Type A for $g$, then, since $\frac{\beta_g}{\gamma_g} = \frac{t}{s}$, we have $y_3 = n_3 s$ and $z_3 = n_3 t$ for some integer $n_3$, and we see that we must have $P(f(n_3)) = P(f(n_2))$, so that we can apply Proposition 3.1 to see that $(a,b,c,x_2,y_2,z_2, x_3, y_3, z_3)$ is in the infinite family (i), contradicting $2 \nmid g > 1$.    
   
So $(x_3,y_3,z_3)$ must be Type B for $g$.  In the part of this proof preceding (4.9) we showed that when 
$(x_2,y_2,z_2)$ is Type A for $g$ we must have Case 1, Case 2, or Case 3; using the same argument under the assumption that $(x_3,y_3,z_3)$ is Type B for $g$, we can show that, noting that now we can assume $2 \nmid b_1$, we must have $b_1 = 1$.  But we have shown that $(a,b,c,x_1,y_1,z_1,x_2,y_2,z_2)$ is in the infinite family (iv) with $d = b_1^{y_1} \ne 1$ by Lemma 8, giving a contradiction.    

So we see that if we have Case 1, then there are exactly two solutions $(x_1,y_1,z_1)$ and $(x_2,y_2,z_2)$.  Thus for this case $N(a,b,c)=2$ and $(a,b,c, x_1,y_1,z_1,x_2,y_2,z_2)$ is in the infinite family (iv) with $iu>1$.  

For Case 2, we have 
$$ a = 2 g^{\alpha_g}, b = g^{\beta_g} b_1, c = g^{\gamma_g} c_1. \eqno{(4.14)} $$
Again letting $ d = b_1^{y_1}$, $b_0 = b^{y_1}$,  $c_0 = c^{z_1}$, we have (4.10).  
$$ b_0 =  g^{\beta_g y_1} d, c_0 = g^{\gamma_g z_1} (d + 2).  \eqno{(4.15)}$$
Note that Case 2 requires $x_1=1$ (by (4.2) and the definition of Case 2). Since the solution $(x_1,y_1, z_1)$ is Type O for $g$, we have $\beta_g y_1 =  \gamma_g z_1  = \alpha_g$.  So we have  
$$ b_0 = g^{\alpha_g} d, c_0 = g^{\alpha_g} (d+2).  \eqno{(4.16)}$$
Letting $k = n_2 = \frac{n_2}{n_1} = \frac{y_2}{y_1} = \frac{z_2}{z_1}$, we see that the solutions $(x_1,y_1,z_1)$ and $(x_2,y_2,z_2)$ for the triple $(a,b,c)$ as in (4.14) correspond to the solutions $(1, 1, 1)$ and $(x_2, k, k)$ for the triple $(a, b_0, c_0)$.  Since $g^{\alpha_g x_2} \parallel a^{x_2} = c_0^k - b_0^k$, from (4.16) we find that $\alpha_g x_2 = k \alpha_g + w_g$ where $g^{w_g} \parallel (d + 2)^k - d^k$ with $w_g > 0$ as in Case 1. So 
$$ x_2 = k + \frac{w_g}{\alpha_g} > 1.  \eqno{(4.17)} $$
Since $2^{x_2} \parallel a^{x_2} = c_0^k - b_0^k$, from (4.16) we see that by Lemma 7 we must have $x_2 = h+v$ where $2^h \parallel 2 d + 2$ and $2^v \parallel k$, $v>0$.  Combining this with (4.17) we obtain
$$ k - v = h - \frac{w_g}{\alpha_g}. $$
Now we see that $(a, b_0, c_0,1, 1, 1, k + ({w_g}/{\alpha_g}), k,k)$ (and therefore also $(a,b,c, x_1,y_1,z_1,x_2,y_2,z_2)$) is in the infinite family (iv) in the Introduction when $i = 1$, $j = \alpha_g$, $u = 1$, and $w = w_g$.  To complete the treatment of Case 2, we need to show that $N(a,b,c)=2$.  It suffices to point out that $(a,b,c)$ gives no third solution $(x_3,y_3,z_3)$ by the same argument used to show there is no third solution $(x_3,y_3,z_3)$ for Case 1.   

So we see that if we have Case 2, then there are exactly two solutions $(x_1,y_1,z_1)$ and $(x_2,y_2,z_2)$.  Thus for this case $N(a,b,c)=2$ and $(a,b,c,x_1,y_1,z_1,x_2,y_2,z_2)$ is in the infinite family (iv) with $iu=1$. 

For Case 3 we have 
$$ a = g^{\alpha_g}, b = g^{\beta_g} b_1, c = g^{\gamma_g} c_1. \eqno{(4.18)} $$
Again letting $ d = b_1^{y_1}$, $b_0 = b^{y_1}$,  $c_0 = c^{z_1}$, we have 
$$ c_1^{z_1} - b_1^{y_1} = c_1^{n_1 t} - b_1^{n_1 s} = f(n_1) = 1. \eqno{(4.19)}$$  
So we have 
$$ b_0 = g^{\beta_g y_1} d, c_0 = g^{\gamma_g z_1} (d + 1).  \eqno{(4.20)}$$
Since the solution $(x_1,y_1, z_1)$ is Type O for $g$, we have $\beta_g y_1 =  \gamma_g z_1  = \alpha_g x_1$.  So we have  
$$ b_0 = g^{\alpha_g x_1} d, c_0 = g^{\alpha_g x_1} (d+1).  \eqno{(4.21)}$$
Again letting $k = n_2$, we find that the solutions $(x_1,y_1,z_1)$ and $(x_2,y_2,z_2)$ for the triple $(a,b,c)$ as in (4.18) correspond to the solutions $(x_1, 1, 1)$ and $(x_2, k, k)$ for the triple $(a, b_0, c_0)$.  Proceeding as in Cases 1 and 2, we find that $\alpha_g x_2 = k \alpha_g x_1 + w_g$ where $g^{w_g} = (d + 1)^k - d^k$, $w_g>0$ . So 
$$ x_2 = k x_1 + \frac{w_g}{\alpha_g}.  \eqno{(4.22)} $$
Now we see that $(a, b_0, c_0, x_1, 1, 1, k x_1 + ({w_g}/{\alpha_g}), k,k)$  (and therefore also $(a,b,c, x_1,y_1,z_1, x_2,y_2,z_2)$) is in the infinite family (iii) in the Introduction, with $a = g^j$, $j = \alpha_g$, $u = x_1$, and $w = w_g$.  (Note that $2 \nmid g$ since $\gcd(g,d(d+1))=1$.)  It remains to show $N(a,b,c)=2$ for Case 3.  We will show that any third solution $(x_3,y_3, z_3)$ must correspond to $(x_2,y_2,z_2)$.  

There can be no third solution which is Type A, C, or O for $g$ by the same argument used for Cases 1 and 2.  And the argument used in the first part of this proof to show that we must have Case 1, 2, or 3 can be used to show that, if there exists a solution $(x_3, y_3,z_3)$ which is Type B for $g$, then $b_1 \le 2$.  If $b_1 = 2$, then, since $a_1=1$, we can reverse the roles of $a$ and $b$ in the argument used for Case 2 to see that $(b,a,c,y_1,x_1,z_1,y_2,x_2,z_2)$ is in the infinite family (iv) with $d=1$, contradicting Lemma 8.     

So if there is a third solution we must have $a_1 = b_1 = 1$.  

From the solution $(x_1,y_1,z_1)$ which is Type O for $g$ we derive $c_1 = 2$.  

From the solution $(x_2,y_2,z_2)$, which is Type A for $g$, we derive
$$ g^{ \alpha_g x_2 - \gamma_g z_2} + 1 = 2^{z_2}. \eqno{(4.23)}$$
It is a familiar elementary result that (4.23) requires $\alpha_g x_2 - \gamma_g z_2 = 1$ and that there are no values of $z$ possible in a solution which is Type A for $g$ other than $z_2$.  

If there is a third solution $(x_3,y_3,z_3)$, it must be Type B for $g$, so that 
$$ 1 + g^{\beta_g y_3 - \gamma_g z_3} = 2^{z_3}.  \eqno{(4.24)}$$
We have $\beta_g y_3 - \gamma_g z_3 = 1$ so that $z_2 = z_3$ is the only possible choice for $z$ in any solution to (1.1) other than $(x_1,y_1,z_1)$.  

Multiplying both sides of (4.23) and (4.24) by $g^{\gamma_g z_2}$ we find that we can apply Criterion 1 to see that $(x_2,y_2,z_2)$ is considered the same as $(x_3,y_3,z_3)$.    

Thus we see that if we have Case 3, then $N(a,b,c)=2$.  Recall also (from the paragraph following (4.22)) that $(a,b,c,x_1,y_1,z_1,x_2,y_2,z_2)$ is in the infinite family (iii). 

Thus we find that in all three cases $N(a,b,c)=2$ with  $(a,b,c,x_1,y_1,z_1,x_2,y_2,z_2)$ in either the infinite family (iii) or the infinite family (iv).  

This completes the proof of Proposition 4.1.  
\end{proof}

Let $S_{aa}$ be the set of all triples $(a,b,c)$ such that $\gcd(a,b)>1$ and (1.1) has two distinct solutions $(x,y,z)$ both of which are Type A for some prime in $Q$.

\begin{Proposition}  
For $(a,b,c) \in S_{aa}$ we have $N(a,b,c)=2$, and (1.1) has two solutions $(x_1,y_1,z_1)$ and$(x_2,y_2,z_2)$ with $(a,b,c, \allowbreak  x_1,y_1,z_1,\allowbreak  x_2,y_2,z_2)$ in one of the infinite families (i), (iii), or (iv). 
\end{Proposition}

\begin{proof} 
For any $(a,b,c)$ with $\gcd(a,b)>1$, (1.1) cannot have two distinct solutions $(x,y,z)$ both of Type A for a given prime which correspond to each other (recall $z_i$ increases with $n_i$ in (3.2)).  So $(a,b,c) \in S_{aa}$ implies $N(a,b,c)>1$.  

Let $J$ be the set of all $(a,b,c) \in S_{aa}$ such that there exists a solution which is Type O for some prime in $Q$, and let $K$ be the set of all $(a,b,c) \in S_{aa}$ such that there are no solutions which are Type O for any prime in $Q$.  

If $(a,b,c) \in J$, then, by Proposition 4.1, $N(a,b,c)=2$ and the equation (1.1) has two solutions $(x_1,y_1,z_1)$ and $(x_2,y_2,z_2)$ with  $(a,b,c,x_1,y_1,z_1,x_2,y_2,z_2)$ in one of the infinite families (iii) or (iv).  

Now suppose $(a,b,c) \in K$.  By the Corollary to Proposition 3.1, $N(a,b,c)=2$ and (1.1) has two solutions $(x_1,y_1,z_1)$ and $(x_2,y_2,z_2)$ with $(a,b,c,x_1,y_1,z_1,x_2,y_2,z_2)$ in the infinite family (i).  

So, since $S_{aa} = J \cup K$, Proposition 4.2 follows.    
\end{proof}

 (Although not needed for our purpose, it is easily seen that when $(a,b,c) \in S_{aa}$ and $(a,b,c, \allowbreak  x_1,y_1,z_1,\allowbreak  x_2,y_2,z_2)$ is in (iv) we must have (iv) with $iu>1$ in order to have $2 \in Q$ so that the two solutions in (iv) are both Type A for 2.  Also, it is not hard to show that when $(a,b,c) \in S_{aa}$, $(a,b,c, \allowbreak  x_1,y_1,z_1,\allowbreak  x_2,y_2,z_2)$ cannot be in the infinite family (iii), but this is not needed for our purpose.)    

\bigskip 

Let $S_{bb}$ be the set of all triples $(a,b,c)$ such that $\gcd(a,b)>1$ and (1.1) has two solutions both of which are Type B for some prime in $Q$.  By the symmetry of $a$ and $b$ (using Definitions 2 and 3), Proposition 4.2 holds with $S_{bb}$ replacing $S_{aa}$.

\begin{Theorem4_3}  
Let $S_i = S_c \cup S_o \cup S_{aa} \cup S_{bb}$.  If $(a,b,c)$ in $S_{i}$, then $N(a,b,c)=2$, and (1.1) has two solutions $(x_1,y_1,z_1)$ and $(x_2,y_2,z_2)$ with $(a,b,c, \allowbreak x_1,y_1,z_1, \allowbreak  x_2, y_2, z_2)$ in one of the infinite families (i), (ii), (iii), or (iv).   
\end{Theorem4_3}

\begin{proof} 
This is an immediate consequence of Propositions 3.4, 4.1, and 4.2.
\end{proof}

Theorem 1.1 follows directly from Theorem 4.3:  

\begin{proof}[Proof of Theorem 1.1]
Eliminate from consideration the cases $(\{a,b\},c)= (\{ 3,5\},2)$ and $(a,b,c) = \allowbreak (2^u, \allowbreak  2^v, \allowbreak 2^w)$, $u, v, w \in \intZ^+$.

For $(a,b,c)$ such that $\gcd(a,b)=1$, Theorem 1 of \cite{MP} shows that $N(a,b,c) \le 2$.  

For $(a,b,c)$ such that $\gcd(a,b)>1$, if we have $N(a,b,c)>1$ and $(a,b,c) \not\in S_i$, then (1.1) has no solutions of Type C or O for any prime in $Q$, and at most one solution of Type A and at most one solution of Type B for any prime in $Q$, so that there are not more than two solutions $(x,y,z)$ to (1.1).  Now Theorem 1.1 is an immediate consequence of Theorem 4.3. 
\end{proof} 

{\bf Comment 4.4 } \quad 
Notice that proving Theorem 1.1 does not require using the statement concerning infinite families in Theorem 4.3 or the similar statements in Propositions 3.1, 3.4, 4.1, and 4.2; only the result $N(a,b,c) = 2$ is needed.  Also, the use of infinite families is not needed for proving $N(a,b,c)=2$  in Propositions 3.1, 3.4, 4.1, and 4.2 except in handling the three special cases in the proof of Proposition 4.1, where using infinite families is convenient but not necessary.  Our primary purpose in considering infinite families is to obtain results concerning anomalous $(a,b,c,x_1,y_1,z_1,x_2,y_2,z_2)$ (as in Theorem 4.5 below) and to obtain results on cases with more than two solutions $(x,y,z)$ to (1.1) when Criterion 1 is not used (as in Theorem 4.7).  The proof of Theorem 4.5 will use the statement concerning infinite families in Theorem 4.3.  Theorem 4.7 will show that the solutions $(x_1,y_1,z_1)$ and $(x_2,y_2, z_2)$ in the statement of Theorem 4.3 (and in the similar statements of Propositions 3.1, 3.4, 4.1, and 4.2)  are the only solutions $(x,y,z)$ to (1.1) for the $(a,b,c)$ in question even when Criterion 1 is not used, except for a few specifically designated $(a,b,c)$ with more than two solutions $(x,y,z)$; for $(a,b,c)$ with more than two solutions $(x,y,z)$ there exists a member $m$ of the set $F$ for a certain infinite family such that any solution to (1.1) for this $(a,b,c)$ must correspond to one of the two solutions $(x, y, z)$ given by $m$.  

\quad 

For Theorem 4.5 and Observation 4.6 which follow we assume $(a,b,c) \ne (2^u, 2^v,2^w)$, $u,v,w \in \intZ^+$.

\begin{Theorem4_5}  
If, for some $(a,b,c)$ with $\gcd(a,b)>1$, (1.1) has two solutions $(x_1,y_1,z_1)$ and $(x_2,y_2, z_2)$ which do not correspond to each other and $(a,b,c, \allowbreak   x_1,y_1,z_1,\allowbreak  x_2,y_2,z_2)$ is not in any of the infinite families (i), (ii), (iii), or (iv), then 
$$ a = g^{\alpha_g} a_1, b = g^{\beta_g} b_1, c = g^{\gamma_g} c_1  \eqno{(4.25)}$$
for some positive integers $g$, $\alpha_g$, $\beta_g$, $\gamma_g$, $a_1$, $b_1$, $c_1$ with $\gcd(g, a_1 b_1 c_1) = 1$ and $\gcd(a_1, b_1) = 1$, and one of the two solutions is Type A for $g$ and the other solution is Type B for $g$, with no further solutions $(x,y,z)$ (even if Criterion 1 is not used).  
\end{Theorem4_5}

\begin{proof}
Assume that, for some $(a,b,c)$ with $\gcd(a,b)>1$, (1.1) has two solutions $(x_1,y_1,z_1)$ and $(x_2,y_2, z_2)$ which do not correspond to each other   and $(a,b,c, \allowbreak   x_1,y_1,z_1,\allowbreak  x_2,y_2,z_2)$ is not in any of the infinite families (i), (ii), (iii), or (iv).  By Theorem 1.1, $N(a,b,c)=2$ and any further solution $(x,y,z)$ must correspond to one of $(x_1,y_1,z_1)$,  $(x_2,y_2,z_2)$, so that no two solutions for this $(a,b,c)$ can correspond to the two solutions given by a member of $F$ for any infinite family (recall the Comment following Definition 2 in the Introduction), so that $(a,b,c) \not\in S_i$ (by Theorem 4.3).   So, recalling the last paragraph in the proof of Theorem 1.1 (immediately following Theorem 4.3), we see that $(x_1,y_1,z_1)$ and $(x_2, y_2, z_2)$ are the only solutions $(x,y,z)$ to (1.1) (even if Criterion 1 is not used), and, for any prime in $Q$, one of $(x_1,y_1,z_1)$ and $(x_2,y_2,z_2)$ is Type A and the other is Type B.  (Note that neither $(x_1,y_1,z_1)$ nor $(x_2,y_2,z_2)$ corresponds to a distinct third solution $(x_3,y_3,z_3)$ since $(a,b,c) \not\in S_{aa} \cup S_{bb}$).     

Suppose one of these solutions, say $(x_1,y_1,z_1)$, is of Type A for some prime $p \in Q$ and of Type B for some prime $q \in Q$.  Then, by the Corollary to Lemma 2, $(x_2,y_2,z_2)$ must also be Type A for $p$, contradicting the previous paragraph (which showed that $(x_2,y_2,z_2)$ must be Type B for $p$).  So we see that one of the solutions, say $(x_1,y_1,z_1)$, must be Type A for every prime in $Q$ and the other solution $(x_2, y_2,z_2)$ must be Type B for every prime in $Q$.    
Let $q_1$ and $q_2$ be any two primes in $Q$.  Then $\beta_{q_1} y_1 = \gamma_{q_1} z_1$, $\beta_{q_2} y_1 = \gamma_{q_2} z_1$, $\alpha_{q_1} x_2 = \gamma_{q_1} z_2$, $\alpha_{q_2} x_2 = \gamma_{q_2} z_2$, from which we derive (3.34), so that we can use Proposition 3.5 to obtain (4.25).  The solution $(x_1,y_1,z_1)$ is Type A for $g$, and the solution $(x_2,y_2,z_2)$ is Type B for $g$.  
\end{proof}

Theorem 4.5 will be used to establish a method for searching for anomalous $(a,b,c, x_1,y_1,z_1,x_2,y_2,z_2)$ in Section 5 which follows.  (Note that in the statement of Theorem 4.5 $(a,b,c, x_1,y_1,z_1,x_2,y_2,z_2)$ is an anomalous nine-tuple.)

Let $I$ be the set of all $(a,b,c)$ with $\gcd(a,b)>1$ such that (1.1) has two solutions $(x_1,y_1,z_1)$ and $(x_2,y_2,z_2)$ with $(a,b,c, x_1,y_1,z_1,x_2,y_2,z_2)$ in one of the four infinite families in the Introduction. 

(Although not needed for what follows, we note that $I = S_i \cup T$ where $T$ is the set of all triples $(a,b,c)$ satisfying either $(\{a,b\}, c) = (\{ 2, 4\}, 2^{2 t + 1} \cdot 3)$ or $(\{4,8\}, 2^{6t+3} \cdot 3)$, where $t \ge 0$.  Since we will not be using this result, we do not give its proof, which uses the method of the proof of Theorem 4.7 below.  Any $(a,b,c) \in T$ gives exactly two solutions $(x_1,y_1,z_1)$ and $(x_2,y_2,z_2)$ to (1.1), one of which is Type A for 2 and the other of which is Type B for 2, with $(a,b,c, x_1, y_1,z_1, x_2,y_2,z_2)$   in the infinite family (i) with $h=2$.)

\begin{Observation4_6}  
The use of Criterion 1 to obtain $N(a,b,c)=2$ is needed only when $(a,b,c)$ gives two solutions $(x_1,y_1,z_1)$ and $(x_2,y_2,z_2)$ to (1.1) with $(a,b,c,x_1,y_1,z_1,x_2,y_2,z_2)$ in either the infinite family (i) with $h=2$ or the infinite family (iii) with $g=2^k -1$, $d =w = j = 1$.  
\end{Observation4_6}  

\begin{proof}
Criterion 1 is relevant only when $a_1=b_1=1$.  Assume $(a,b,c)$ satisfies $\gcd(a,b)>1$, $N(a,b,c)=2$, and $a_1 = b_1 = 1$.  

If $(a,b,c) \not\in I$, then $(a,b,c)$ satisfies the conditions of the statement of Theorem 4.5, so, by Theorem 4.5, $(a,b,c)$ has exactly two solutions $(x,y,z)$, neither of which corresponds to a further distinct solution $(x,y,z)$.  So Criterion 1 does not apply here.   

If $(a,b,c) \in I$, then, since $a_1=b_1=1$, it suffices to determine when any of the infinite families (i), (ii), (iii), (iv) allows $a_1=b_1=1$.  In the infinite family (ii), we have $a_1 = 2$, so we can eliminate (ii) from consideration.  In the infinite family (iv), we cannot have $b_1 = d = 1$ by Lemma 8, so we can eliminate (iv) from consideration.  In the infinite family (i), $b_1 = 1$ holds only if $h = 2$.  Finally, in the infinite family (iii), $b_1 = 1$ requires $d=1$, so that, since $g^w = (d+1)^k - d^k = 2^k - 1$, it is a familiar elementary result that  $g=2^k - 1$ and $w=1$, so that, since $x_2 = ku + (w/j)$, we must have $j=1$.  
\end{proof}

Using Observation 4.6 we can consider (1.1) without using Criterion 1 and determine all cases for which (1.1) has more than two solutions $(x,y,z)$.

\begin{Theorem4_7}   
Let $(a,b,c)$ be a triple giving more than two solutions $(x,y,z)$ to (1.1).  

If $\gcd(a,b)=1$, then $(\{ a,b\},c) = (\{3,5 \}, 2)$, and (1.1) has exactly three solutions:  $3+5=2^3$, $3^3+5=2^5$, $3+5^3 = 2^7$.  

If $\gcd(a,b)>1$, then we must have one of the following:

1.)  $(a,b,c) = (2,2,2^\gamma \cdot 3)$, $\gamma \in \intZ^+$, for which (1.1) has exactly four solutions $(x,y,z)$: $(\gamma+1, \gamma, 1)$, $(\gamma, \gamma+1, 1)$, $(2 \gamma +3, 2 \gamma, 2)$, $(2 \gamma, 2 \gamma +3, 2)$.  

2.)  $(a,b,c)=(2, 8, 2^{3t} 3)$, $t \in \intZ^+$, for which (1.1) has exactly three solutions $(x,y,z)$: $(3t+1, t, 1)$, $(6t+3, 2t, 2)$, $(6t, 2t+1, 2)$

or

$(a,b,c) = (8,2,2^{3t} 3)$, $t \in \intZ^+$, for which (1.1) has exactly three solutions $(x,y,z)$: $(t, 3t+1, 1)$, $(2t, 6t+3, 2)$, $(2t+1, 6t, 2)$.  

3.)  $(a,b,c) = (2^k - 1, 2^k - 1, 2 (2^k - 1)^\gamma)$, $k, \gamma \in \intZ^+$, for which (1.1) has exactly three solutions $(x,y,z)$: $(\gamma, \gamma, 1)$, $(k \gamma +1, k \gamma, k)$, $(k \gamma, k \gamma +1,k)$.  

4.)  $(a,b,c) = (2^u, 2^v, 2^w)$, $u, v, w \in \intZ^+$, $\gcd(uv, w)=1$, 
for which (1.1) has an infinite number of solutions $(x,y,z) = (\frac{t v}{g}, \frac{tu}{g}, \frac{t L + 1}{w})$ where $\gcd(u,v)=g$, $\lcm(u,v)=L$, and $t$ is a positive integer such that $tL \equiv -1 \bmod w$.  
\end{Theorem4_7}

\begin{proof} 
If $\gcd(a,b)=1$, then the result that $(\{a,b\},c) = (\{3,5\},2)$ is given by \cite{MP}, and the result that there are exactly three solutions in this case is given by \cite{MP} as well as earlier papers cited in \cite{MP}.  

Now assume $(a,b,c)$ with $\gcd(a,b)>1$ gives more than two solutions $(x,y,z)$ to (1.1).  If $a$, $b$, and $c$ are all powers of 2 and (1.1) has a solution, then we must have {\it 4.)}; assume $a$, $b$, and $c$ are not all powers of 2, so that we are eliminating {\it 4.)} from consideration. Since any solution to (1.1) has at most one solution distinct from it which corresponds to it, we must have $N(a,b,c)>1$.  By Theorem 1.1, $N(a,b,c)=2$.  

By Observation 4.6, we must have one of the following two cases. 

{\it Case 1}: $(a,b,c)$ gives two solutions $(x_1,y_1,z_1)$, $(x_2,y_2,z_2)$ to (1.1) such that $(a,b,c,x_1,y_1,z_1,x_2,y_2,z_2)$ is in the infinite family (i) with $h=2$, and there is at least one further solution $(x_3,y_3,z_3)$ where the solution $(x_3,y_3,z_3)$ corresponds to one of $(x_1,y_1,z_1)$ and $(x_2,y_2,z_2)$ (recall Criterion 1, noting that $N(a,b,c)=2$ and the solutions $(x_1,y_1,z_1)$ and $(x_2,y_2,z_2)$ do not correspond to each other since, in any infinite family, the two solutions $(x,y,z)$ do not correspond to each other).  

For this case, we have $(a,b,c) = (2^\alpha, 2^\beta, 2^\gamma \cdot 3)$ since, by Definition 3, the solutions $(x_1,y_1,z_1)$ and $(x_2,y_2,z_2)$ must correspond to the two solutions given by some member of $F$ for the infinite family (i) (note that this requires $\min(z_1,z_2) =1$).  Letting $m$ be this member of $F$ we have $m = (2, 2^u, 2^u \cdot 3, u+1, 1, 1, 2u+3, 2, 2)$ for some fixed positive integer $u$.  By Definition 3, $\gamma = \gamma \min(z_1, z_2) = u \min(1,2) = u$, so that 
$$m = (2,2^\gamma, 2^\gamma \cdot 3, \gamma+1, 1, 1, 2 \gamma+3, 2,2). $$ 
Let $(x_1,y_1,z_1)$ be the solution to (1.1) which corresponds to the solution $(\gamma+1, 1, 1)$ in $m$, and let $(x_2,y_2,z_2)$ be the solution which corresponds to the solution $(2 \gamma + 3, 2,2)$ in $m$.  Then we have
$$ \{ \alpha x_1, \beta y_1 \} = \{ \gamma+1, \gamma\}  \eqno{(4.26)}$$
and 
$$ \{ \alpha x_2, \beta y_2 \} = \{ 2 \gamma + 3, 2 \gamma \}.  \eqno{(4.27)}$$
Let $(x_3,y_3,z_3)$ be a further solution.  If $(x_3,y_3,z_3)$ corresponds to $(x_1,y_1,z_1)$, then, if $\alpha x_1 = \gamma +1$ (respectively, $\gamma$), we must have $\alpha x_3 = \gamma$ (respectively, $\gamma+1$).  By (4.26) we see that this requires $\alpha \mid \gamma+1$, $\alpha \mid  \gamma$, $\beta  \mid  \gamma+1$, $\beta \mid \gamma$, so that $\alpha = \beta = 1$.  

If $(x_3,y_3,z_3)$ corresponds to $(x_2,y_2,z_2)$, then by (4.27) we must have $\alpha \mid 2 \gamma + 3$, $\alpha \mid 2 \gamma$, $\beta \mid 2 \gamma + 3$, $\beta \mid 2 \gamma$.  This requires $\alpha \in \{ 1, 3\}$, $\beta \in \{ 1,3 \}$.  so we have $(\alpha, \beta ) = (1,1)$, $(1,3)$, or $(3,1)$, noting that $(\alpha, \beta) = (3,3)$ is impossible by (4.26).  

 If $\alpha = \beta = 1$ we have {\it 1.)} in the formulation of Theorem 4.7.  
If $(\alpha, \beta) = (1,3)$ or $(3,1)$, we have {\it 2.)} in the formulation of Theorem 4.7 (note that in this case we must have $3 \mid \gamma$).   

{\it Case 2}:  $(a,b,c)$ gives two solutions $(x_1,y_1,z_1)$, $(x_2,y_2,z_2)$ to (1.1) such that $(a,b,c,x_1,y_1,z_1,x_2,y_2,z_2)$ is in the infinite family (iii) with $g = 2^k - 1$, $d=w=j=1$, and there is at least one further solution $(x_3,y_3,z_3)$ where the solution $(x_3,y_3,z_3)$ corresponds to one of $(x_1,y_1,z_1)$ and $(x_2,y_2,z_2)$ as in Case 1. 

For this case we have $(a,b,c)= ((2^k -1)^\alpha, (2^k-1)^\beta, 2 (2^k-1)^\gamma)$ since, by Definition 3, the solutions $(x_1,y_1,z_1)$ and $(x_2,y_2,z_2)$ must correspond to the two solutions given by some member of $F$ for the infinite family (iii) (note that this requires $\min(z_1,z_2) = 1$, and recall the elementary result that $2^k -1$ is never a perfect power).  Letting $m$ be this member of $F$ and proceeding as in Case 1, we see that we must have 
$$m = (2^k-1, (2^k-1)^\gamma, 2 (2^k-1)^\gamma, \gamma, 1, 1, k \gamma + 1, k, k)$$
for some positive integer $k$.  

Let $(x_1,y_1,z_1)$ be the solution to (1.1) which corresponds to the solution $(\gamma,1,1)$ in $m$, and let $(x_2,y_2,z_2)$ be the solution to (1.1) which corresponds to the solution  $(k\gamma+1, k,k)$.  

Let $(x_3,y_3,z_3)$ be a further solution.  If $(x_3,y_3,z_3)$ corresponds to $(x_1,y_1,z_1)$, then both $(x_1,y_1,z_1)$ and $(x_3,y_3,z_3)$ are Type O for $g=2^k-1$, so that $\alpha x_1 = \beta y_1 = \alpha x_3 = \beta y_3$.  This requires $x_1=x_3$, $y_1= y_3$, so that $(x_3,y_3,z_3)$ is not a distinct solution.  

If $(x_3,y_3,z_3)$ corresponds to $(x_2,y_2,z_2)$, then we can proceed as in Case 1 to see that $\alpha \mid k \gamma+1$, $\alpha \mid k \gamma$, $\beta \mid k \gamma +1$, $\beta  \mid  k \gamma$, so that $\alpha = \beta = 1$, which gives {\it 3.)}  in the formulation of Theorem 4.7. 
\end{proof} 

From Theorem 4.7 we can immediately obtain a revised version of Theorem 1.1 in which Criterion 1 is replaced by a more specific restriction:

\begin{Theorem4_8} 
For given integers $a$, $b$, $c$ all greater than one, (1.1) has at most two solutions in positive integers $(x,y,z)$ (where two solutions $(x_1,y_1,z_1)$, $(x_2 , y_2, z_2)$  are considered the same solution if $a=b$ and $\{ x_1,y_1 \} = \{ x_2, y_2\}$), except for the following three cases: 

1.)  $\{ a,b\} = \{ 3,5\}$, $c=2$ which gives the three solutions in Theorem 4.7.  

2.)  $\{ a,b\} = \{ 2,8 \}$, $c=2^{3t} 3$, which gives the three solutions in 2.) of Theorem 4.7.  

3.)  $\{ a,b\} = \{ 2^u, 2^v\}$, $c = 2^w$, where $\gcd(uv,w)=1$, which gives the infinite number of solutions in 4.) of Theorem 4.7.
\end{Theorem4_8}

\section{Cases with exactly two solutions} 

By Theorem 1.1 we have $N(a,b,c) \le 2$ for all $(a,b,c)$ except $(\{a,b\}, c) = (\{3,5\},2)$ or $(\{ 2^u, 2^v\}, 2^w)$ for positive integers $u$, $v$, $w$ (where $\gcd(uv,w)=1$ by {\it 4.)} of Theorem 4.7).  
Let $S_j$ be the set of all triples $(a,b,c)$ such that $\gcd(a,b)>1$, $N(a,b,c)=2$, and there do not exist solutions $(x_1,y_1,z_1)$ and $(x_2,y_2,z_2)$ to (1.1) with $(a,b,c,x_1,y_1,z_1,x_2,y_2,z_2)$ in any of the infinite families (i), (ii), (iii), or (iv).  By Theorem 4.5, for any $(a,b,c) \in S_j$ there are exactly two solutions $(x,y,z)$ one of which is Type A for $g$, the other of which is Type B for $g$, and neither of which has a solution corresponding to it (here $g$ is as in Theorem 4.5).  In this section we consider whether there exist any $(a,b,c) \in S_j$ which are not listed among the ten anomalous cases given in the Introduction.  (Note that, by Theorem 4.5, $S_j$ is the set of $(a,b,c)$ for which anomalous nine-tuples occur.)  

We can assume that for any $(a,b,c) \in S_j$, we have a solution $(x_1,y_1,z_1)$ which is Type A for $g$ and a solution $(x_2,y_2,z_2)$ which is Type B for $g$, where $g$ is as in Theorem 4.5.  Using the notation of Theorem 4.5 (taking $a_1 \le  b_1$ and, for brevity, writing $\alpha$ for $\alpha_g$, $\beta$ for $\beta_g$, and $\gamma$ for $\gamma_g$), from the solutions $(x_1,y_1,z_1)$ and $(x_2,y_1,z_2)$ we derive two equations with relatively prime terms:    
$$ g^{\alpha  x_1 - \gamma  z_1} a_1^{x_1} + b_1^{y_1} = c_1^{z_1} \eqno{(5.1)}$$
and
$$ a_1^{x_2} + g^{\beta  y_2 - \gamma  z_2} b_1^{y_2} = c_1^{z_2}. \eqno{(5.2)}$$
For any $(a,b,c) \in S_j$, the two solutions $(x_1,y_1,z_1)$ and $(x_2,y_2,z_2)$ must be derived from equations with relatively prime terms as in (5.1) and (5.2).  So, in searching for hitherto unknown anomalous solutions we can begin by examining pairs of equations of the form 
$$ g^{w_1} a_1^{x_1} + b_1^{y_1} = c_1^{z_1} \eqno{(5.3)}$$
and 
$$a_1^{x_2} + g^{w_2} b_1^{y_2} = c_1^{z_2} \eqno{(5.4)}$$
where $g$, $a_1$, $b_1$, and $c_1$ are pairwise relatively prime.  
For any pair of equations with relatively prime terms which can be represented as in(5.3) and (5.4), we consider whether there exist $\alpha$, $\beta$, and $\gamma$ which produce an $(a,b,c) \in S_j$.  To do this we construct systems of linear equations in which the unknowns are $\alpha$, $\beta$, and $\gamma$:
$$ y_1 \beta = z_1 \gamma, \eqno{(5.5)}$$
$$ y_2 \beta - w_2 = z_2 \gamma, \eqno{(5.6)}$$
$$ x_1 \alpha - w_1 = z_1 \gamma, \eqno{(5.7)}$$
$$  x_2 \alpha = z_2 \gamma.  \eqno{(5.8)}$$

If $1 < a_1 < b_1$, we have a system of four equations which must have a solution in three positive integer variables $\alpha$, $\beta$, $\gamma$ in order to produce an $(a,b,c) \in S_j$.  

If $1 = a_1 < b_1$, then it suffices to obtain a solution in positive integers $\beta$ and $\gamma$ for the system of two equations (5.5) and (5.6) in order to produce an $(a,b,c) \in S_j$, since then we can take $\alpha = 1$ and let 
$$ x_1 = z_1 \gamma + w_1 \eqno{(5.9)}$$
and 
$$ x_2 = z_2 \gamma.  \eqno{(5.10)}$$

If $1 = a_1 = b_1$, then we have 
$$ g^{w_1} + 1 = c_1^{z_1} \eqno{(5.11)}$$
and 
$$ 1 + g^{w_2} = c_1^{z_2}.  \eqno{(5.12)}$$
If $z_1 = z_2$, then the solutions $(x_1,y_1,z_1)$ and $(x_2,y_2,z_2)$ to (1.1) from which (5.11) and (5.12) are derived must correspond to each other, so, by Observation 4.6, $(a,b,c) \not\in S_j$.  So $z_1 \ne z_2$, in which case, since (5.12) can also be considered an equation derived from a solution of Type A for $g$, we have $(a,b,c) \in S_{aa} \subset S_i$, contradicting the definition of $S_j$.  (Although not needed for our purposes, we note that, when $z_1 \ne z_2$, (5.11) and (5.12) must be $2+1=3$ and $2^3+1=3^2$.)

In conducting a search for possible further anomalous solutions not already known, we can begin by considering ternary equations with relatively prime terms.  A remarkably comprehensive list of such equations was constructed by Matschke \cite{Mat} based on work of K\"anel and Matschke \cite{KM}. This list gives all cases of $A+B=C$ with $\gcd(A,B)=1$ and $\rad(ABC) < 10^7$.  From this list we find pairs of equations satisfying (5.3) and (5.4) and determine whether a triple of positive integers $(\alpha, \beta, \gamma)$ satisfies (5.5) through (5.8).  We thus show that there are no $(a,b,c)$ with $\gcd(a,b)>1$, $N(a,b,c)=2$, and $\rad(abc)<10^7$ other than those already known and listed among the ten anomalous cases and the four infinite families in the Introduction.  

(One can calculate that there are no further pairs (5.3) and (5.4) allowing $\alpha$, $\beta$, $\gamma$ which satisfy (5.5) through (5.8) for $a_1 \le 100$. $g \le 100$,  $b_1 \le 10000$, and the exponents $w_1$, $x_1$, $y_1$, $z_1$, $x_2$, $w_2$, $y_2$, $z_2$ each less than or equal to 10.)   
         
These results complement the bounds given in the Introduction found by an independent computer search.

\section*{Acknowledgments} 

We are grateful to Takafumi Miyazaki for many helpful suggestions and corrections.  This work used the Augie High Performance Computing cluster, funded by award NSF 2018933, at Villanova University.


\begin{thebibliography}{1}


\bibitem{Be}
M. Bennett, On some exponential equations of S. S. Pillai, {\it Canadian J Math.}, {\bf 53} (2001), 897--922.  

\bibitem{BSch} 
F. Beukers and H.P. Schlickewei,
The equation $x + y = 1$ in finitely generated groups,
{\it Acta Arith.},
{\bf 78} (1996),
189---199.


\bibitem{G}
A.~O. {Gel'fond},
\newblock Sur la divisibilit\'e de la diff\'erence des puissance de deux
  nombres entiers par une puissance d'un id\'eal premier,
\newblock {\em Mat. Sb.}, {\bf 7} (49) (1940), 7--25.


\bibitem{HK} 
N. Hirata-Kohno, S-unit equations and integer solutions to exponential Diophantine equations, in {\it Analytic Number Theory and surrounding Areas 2006}, Kyoto RIMS Kokyuroku, 2006, 92--97.


\bibitem{HL1}
Y.-Z. Hu and M.-H. Le,
\newblock An upper bound for the number of solutions of ternary purely
  exponential diophantine equations,
\newblock {\em J. Number Theory}, {\bf 183} (2018), 62--73.

\bibitem{HL2} 
Y.-Z. Hu and M.-H. Le, 
An upper bound for the number of solutions of ternary purely exponential diophantine equations II, 
{\it Publ. Math. Debrecen} {\bf 95} (2019),
335--354.


\bibitem{KM} 
R. von K\"anel and B. Matschke, Solving S-unit, Mordell, Thue, Thue-Mahler and generalized Ramanujan-Nagell equations via Shimura-Taniyama conjecture, Mem. Amer. Math. Soc. {\bf 286} (2023), no.1419, vi+142 pp.


\bibitem{LS}
Maohua Le and Robert Styer, On a conjecture concerning the number of solutions to $a^x+b^y=c^z$, {\it Bulletin of the Australian Mathematical Society}, {\bf 108}  (2023).



\bibitem{LSS}
Maohua Le and Reese Scott and Robert Styer, On a conjecture concerning the number of solutions to $a^x+b^y=c^z$, II. to appear {\it Glasnik Math.}.



\bibitem{M}
K.~Mahler,
\newblock {Zur Approximation algebraischer Zahlen I: \"Uber den gr\"ossten
  Primteiler bin\"arer Formen},
\newblock {\em Math. Ann.}, {\bf 107} (1933), 691--730.

\bibitem{Mat}
B. Matschke, Data attached to the paper \lq Solving S-unit, Mordell, Thue, Thue-Mahler and generalized Ramanujan-Nagell equations via Shimura-Taniyama conjecture' \url{ https://www.math.u-bordeaux.fr/~bmatschke/data/}



\bibitem{MP}
T. Miyazaki and I. Pink, Number of solutions to a special type of unit equation in two variables, {\it Amer. J. Math.} {\bf 146} (2024), no. 2, 295--369.


\bibitem{MP2}
T. Miyazaki and I. Pink, Number of solutions to a special type of unit equation in two variables, II. {\it Research in Number Theory},  {\bf 10} (2024), no. 2, 41 p.


\bibitem{MP3}
T. Miyazaki and I. Pink, Number of solutions to a special type of unit equation in two variables, III. (2024) arXiv:2403.20037.   



\bibitem{Sc}
R. Scott, 
\newblock On the Equations $p^x-b^y = c$ and $a^x+b^y=c^z$, 
\newblock {\it Journal of Number Theory}, {\bf 44}, no. 2 (1993), 153--165. 



\bibitem{ScSt6}
R. Scott and R. Styer. Bennett’s Pillai theorem with fractional bases and negative exponents allowed.
{\it Journal de théorie des nombres de Bordeaux}, {\bf 27} no. 1 (2015), 289--307. 

\bibitem{ScSt6a}     
R. Scott and R. Styer, 
Number of solutions to $a^x + b^y = c^z$, 
{\em Publ. Math. Debrecen}
{\bf 88}
(2016),
 131--138.


\bibitem{W}
A. Wiles. Modular elliptic curves and Fermat’s last theorem. {\it Ann.of Math.},{\bf 141}. (1995) 443--551.




\end{thebibliography}
\end{document}